\documentclass[10pt]{amsart}
\usepackage[T1]{fontenc}
\usepackage{geometry}
\usepackage[latin1] {inputenc}
\usepackage{amsmath}
\usepackage{amsfonts, amssymb, textcomp}
\usepackage[colorlinks=flase, linkcolor=red,urlcolor=green, citecolor=blue]{hyperref}
\usepackage{subeqnarray}
\usepackage{color}

\usepackage{latexsym}
\usepackage{fancyhdr}
\usepackage{longtable}
\usepackage{amsmath, amssymb}
\usepackage{graphicx}

\numberwithin{equation}{section}

\newtheorem{lemma}[subsection]{Lemma}
\newtheorem{theorem}[subsection]{Theorem}
\newtheorem{proposition}[subsection]{Proposition}
\newtheorem{corollary}[subsection]{Corollary}
\newtheorem{definition}[subsection]{Definition}
\newtheorem{example}[subsection]{Example}

\newcommand{\RR}{\mathbb{R}}
\newcommand{\CC}{\mathbb{C}}
\newcommand{\NN}{\mathbb{N}}

\newcommand{\sym}{\operatorname{sym}}

\newcommand{\ev}{\operatorname{ev}}
\newcommand{\ins}{\operatorname{ins}}
\newcommand{\comp}{\operatorname{comp}}

\let\on=\operatorname

\newenvironment{demo}[1]{\par\smallskip\noindent{\bf #1.}}{\par\smallskip}

\title[The convenient setting for ultradifferentiable functions]
{The convenient setting for ultradifferentiable mappings of Beurling- and Roumieu-type defined by a weight matrix}

\author[G.~Schindl]{Gerhard Schindl}

\address{G.~Schindl: Fakult\"at f\"ur Mathematik, Universit\"at Wien,
Oskar-Morgenstern-Platz~1, A-1090 Wien, Austria}
\email{a0304518@unet.univie.ac.at}

\begin{document}
\begin{abstract}
We prove in a uniform way that all ultradifferentiable function classes $\mathcal{E}_{\{\mathcal{M}\}}$ of Roumieu-type and $\mathcal{E}_{(\mathcal{M})}$ of Beurling-type defined in terms of a weight matrix $\mathcal{M}$ admit a convenient setting if $\mathcal{M}$ satisfies some mild regularity conditions. For $\mathcal{C}$ denoting either $\mathcal{E}_{\{\mathcal{M}\}}$ or $\mathcal{E}_{(\mathcal{M})}$ the category $\mathcal{C}$ is cartesian closed, i.e. $\mathcal{C}(E\times F,G)\cong\mathcal{C}(E,\mathcal{C}(F,G))$ for $E,F,G$ convenient vector spaces. As special cases one obtains the classes $\mathcal{E}_{\{M\}}$ and $\mathcal{E}_{(M)}$ respectively $\mathcal{E}_{\{\omega\}}$ and $\mathcal{E}_{(\omega)}$ defined by a weight sequence $M$ respectively a weight function $\omega$.
\end{abstract}

\thanks{GS was supported by FWF-Project P~23028-N13 and FWF-Project P~26735-N25}
\date{\today}

\maketitle
\section{Introduction}
Spaces of ultradifferentiable functions are subclasses of smooth functions with certain growth conditions on all their derivatives. In the literature two different approaches are considered, either using a weight sequence $M=(M_k)_k$ or using a weight function $\omega$. For compact $K$ the set
$$\left\{\frac{f^{(k)}(x)}{h^kM_k}: x\in K, k\in\NN\right\}\hspace{30pt}\text{respectively}\hspace{30pt}\left\{\frac{f^{(k)}(x)}{\exp(1/l\varphi^{*}_{\omega}(lk))}: x\in K, k\in\NN\right\}$$
should be bounded, where the positive real number $h$ respectively $l$ is subject to either a universal or an existential quantifier and $\varphi^{*}_{\omega}$ denotes the Young-conjugate of $\varphi_{\omega}=\omega\circ\exp$. In the case of a universal quantifier we call the class of {\itshape Beurling-type,} denoted by $\mathcal{E}_{(M)}$ or $\mathcal{E}_{(\omega)}$, in the case of an existential quantifier we call the class of {\itshape Roumieu-type,} denoted by $\mathcal{E}_{\{M\}}$ or $\mathcal{E}_{\{\omega\}}$. We write $\mathcal{E}_{[\star]}$ if either $\mathcal{E}_{\{\star\}}$ or $\mathcal{E}_{(\star)}$ is considered.

That a class of mappings $\mathcal{C}$ admits a {\itshape convenient setting} means that one can extend the class to admissible infinite dimensional vector spaces $E,F,G$ such that $\mathcal{C}(E,F)$ is again admissible and the spaces $\mathcal{C}(E\times F,G)$ and $\mathcal{C}(E,\mathcal{C}(F,G))$ are canonically $\mathcal{C}$-diffeomorphic. This important property is called the {\itshape exponential law.}\vspace{6pt}

We recall now some facts, see \cite{KM97} or the appendix in \cite{KMRc} for a short overview. The class $\mathcal{E}$ of all smooth functions admits a convenient setting and for this approach one can test smoothness along $\mathcal{E}$-curves. The class $\mathcal{C}^{\omega}$ of all real-analytic mappings also admits a convenient setting. A mapping is $\mathcal{C}^{\omega}$ if and only if it is $\mathcal{E}$ and in addition it is weakly $\mathcal{C}^{\omega}$ along (weakly) $\mathcal{C}^{\omega}$-curves, i.e. curves whose compositions with any bounded linear functional are $\mathcal{C}^{\omega}$. It actually suffices to test along affine lines.\vspace{6pt}

In \cite{KMRc}, \cite{KMRq} and finally in \cite{KMRplot} A. Kriegl, P.W. Michor and A. Rainer were able to develop the convenient setting for all reasonable classes $\mathcal{E}_{(M)}$ and $\mathcal{E}_{\{M\}}$. In the first step in \cite{KMRc} they introduced the convenient setting for $\mathcal{E}_{\{M\}}$ by testing with $\mathcal{E}_{\{M\}}$-curves for {\itshape non-quasianalytic, strongly log-convex} weight sequences $M$ of {\itshape moderate growth.} A function is $\mathcal{E}_{\{M\}}$ if and only if it is $\mathcal{E}_{\{M\}}$ along all $\mathcal{E}_{\{M\}}$-curves. It was shown that moderate growth is really necessary for the exponential law and non-quasianalyticity is needed for the existence of $\mathcal{E}_{\{M\}}$-partitions of unity.

Then, in \cite{KMRq}, they succeeded to introduce the convenient setting for some quasianalytic classes $\mathcal{E}_{\{M\}}$. In this case $M$ has to satisfy again strong log-convexity, moderate growth and be such that $\mathcal{E}_{\{M\}}$ can be represented as the intersection of all larger non-quasianalytic classes $\mathcal{E}_{\{L\}}$ with strongly log-convex $L$. A mapping is $\mathcal{E}_{\{M\}}$ if and only if it is $\mathcal{E}_{\{L\}}$ along each $\mathcal{E}_{\{L\}}$-curve for each $L\ge M$ which is strongly log-convex and non-quasianalytic. A family of explicit examples $\mathcal{E}_{\{M\}}$ satisfying the requested assumptions was constructed, but the approach does not cover the real analytic case $\mathcal{C}^{\omega}$ and thus was not completely satisfactory.\vspace{6pt}

Finally, in \cite{KMRplot}, it was shown that all classes $\mathcal{E}_{\{M\}}$ and $\mathcal{E}_{(M)}$ such that $M$ is {\itshape strongly log-convex} and has {\itshape moderate growth} admit a convenient setting, no matter if $M$ is quasianalytic or not. Instead of testing along curves the mappings are tested along Banach plots, i.e. mappings of the respective weak class defined in open subsets of Banach spaces. A smooth mapping between convenient vector spaces is $\mathcal{E}_{[M]}$ if it maps $\mathcal{E}_{[M]}$-Banach-plots to $\mathcal{E}_{[M]}$-Banach-plots.\vspace{6pt}

The aim of this work is to generalize the results of \cite{KMRplot} to classes $\mathcal{E}_{[\mathcal{M}]}$ defined by (one-parameter) weight matrices $\mathcal{M}:=\{M^x: x\in\RR_{>0}\}$. In \cite{compositionpaper} the classes $\mathcal{E}_{[M]}$ and $\mathcal{E}_{[\omega]}$ were identified as particular cases of $\mathcal{E}_{[\mathcal{M}]}$. So using this new approach one is able to transfer results from one setting into the other one. Moreover one is able to prove results for $\mathcal{E}_{[M]}$ and $\mathcal{E}_{[\omega]}$ simultaneously and no longer two separate proofs are necessary. We have also shown that there are classes $\mathcal{E}_{[\mathcal{M}]}$ which cannot be described by a single $M$ or $\omega$, e.g. the class defined by the {\itshape Gevrey-matrix} $\mathcal{G}:=\{(p!^{s+1})_{p\in\NN}: s>0\}$. To transfer the proofs of \cite{KMRplot} we will assume for $\mathcal{M}$ among mild basic properties the so-called generalized {\itshape Faà-di-Bruno-property} $(\mathcal{M}_{[\text{FdB}]})$ and the {\itshape moderate growth condition} $(\mathcal{M}_{[\text{mg}]})$.\vspace{6pt}

After introducing the basic notation and definitions we recall the setting of Whitney jets between Banach spaces. We introduce classes of ultradifferentiable functions defined by weight matrices, first between Banach spaces and then between convenient vector spaces. This will be done in section \ref{section2}. In section \ref{section3} we are going to prove the most important and new tools in this work. We will develop projective descriptions for the classes $\mathcal{E}_{[\mathcal{M}]}$ in order to get rid of both existence quantifiers in the Roumieu-case (if $\mathcal{M}=\{M\}$ only one occurs). For this we have to use diagonal techniques and to introduce several families of sequences of positive real numbers to generalize the results of \cite{KMRplot}. These projective representations are needed in section \ref{section4} for the proof of Theorem \ref{Banach-plots-category} to show that $\mathcal{E}_{[\mathcal{M}]}$ is a category and for cartesian closedness Theorem \ref{exponentiallaw} in section \ref{section5}.\vspace{6pt}

Finally in section \ref{section6} we summarize some special cases. In \ref{section63} we revisit weight matrices as defined by Beaugendre in \cite{Beaugendre01} and Schmets and Valdivia in \cite{intersectionpaperextension}. Put $\mathcal{M}^{\Phi}:=\{(p!m^{\Phi}_{ap})_{p\in\NN}: a>0\}$, where $\Phi:[0,+\infty)\rightarrow\RR$ is a convex and increasing function with $\lim_{t\rightarrow\infty}\frac{\Phi(t)}{t}=+\infty$, $\Phi(0)=0$. In the literature only the Beurling-type-class was studied. We will see that the results in this work can also be applied to such classes.\vspace{6pt}

Note that if $\mathcal{M}=\{M\}$ then the Faà-di-Bruno-property for $M$ is sufficient to show closedness under composition and is sufficient for the proofs in this work. But it is really weaker than strong log-convexity as assumed always in the previous papers and proofs of Kriegl, Michor, Rainer, see \cite[3.3.]{compositionpaper} for an explicit (counter)-example. So our results are slightly more general than those of \cite{KMRplot} even in the single weight sequence case. In Lemma \ref{moderategrowthcounterexample} we will show that $(\mathcal{M}_{\{\text{mg}\}})$ is necessary for cartesian closedness of $\mathcal{E}_{\{\mathcal{M}\}}$ and in Example \ref{Matsumoto-counterexample} we will point out that there exist weight matrices $\mathcal{M}$ such that no $M^x\in\mathcal{M}$ has moderate growth but nevertheless $(\mathcal{M}_{\{\text{mg}\}})$ is valid. In particular this holds if the matrix is associated to a weight function $\omega$ and such that $\mathcal{E}_{[\omega]}=\mathcal{E}_{[M]}$ does not hold, see \cite{BonetMeiseMelikhov07} and \cite{compositionpaper}.\vspace{6pt}

This paper contains some of the main results of the authors PhD-Thesis, see \cite{dissertation}. The author thanks his advisor A. Kriegl, P.W. Michor and A. Rainer for the supervision and their helpful ideas.

\subsection{Basic notation}
We denote by $\mathcal{C}$ the class of all continuous, by $\mathcal{E}$ the class of smooth functions and $\mathcal{C}^{\omega}$ is the class of all real analytic functions. We will write $\NN_{>0}=\{1,2,\dots\}$, $\NN=\NN_{>0}\cup\{0\}$ and put $\RR_{>0}:=\{x\in\RR: x>0\}$. For $\alpha=(\alpha_1,\dots,\alpha_n)\in\NN^n$ we use the usual multi-index notation, write $\alpha!:=\alpha_1!\dots\alpha_n!$, $|\alpha|:=\alpha_1+\dots+\alpha_n$ and for $x=(x_1,\dots,x_n)\in\RR^n$ we set $x^{\alpha}=x_1^{\alpha_1}\cdots x_n^{\alpha_n}$. We also put $\partial^{\alpha}=\partial_1^{\alpha_1}\cdots\partial_n^{\alpha_n}$ and denote by $f^{(k)}$ the $k$-th order {\itshape Fréchet derivative} of $f$. Iterated uni-directional derivatives are defined by
$d^k_vf(x):=\left(\frac{d}{dt}\right)^k f(x+tv)|_{t=0}$.

Let $E_1,\dots,E_k$ and $F$ be topological vector spaces, then $L(E_1,\dots,E_k,F)$ is the space of all bounded $k$-linear mappings $E_1\times\dots\times E_k\rightarrow F$. If $E=E_i$ for $i=1,\dots,k$, then we write $L^k(E,F)$. $L^k_{\sym}(E,F)$ is the space of all symmetric $k$-linear bounded mappings $\underbrace{E\times\dots\times E}_{k-\text{times}}\rightarrow F$, so $f^{(k)}:U\rightarrow L^k_{\sym}(E,F)$. $E^{\ast}$ denotes the space of all continuous linear functionals on $E$, $E^{'}$ the space of all bounded linear functionals. If $B\subseteq E$ is closed absolutely convex bounded, then $E_B$ denotes the space generated by $B$ with the {\itshape Minkowski-functional} $\|\cdot\|_B$.

Let $E$ be a locally convex vector space, then the $c^{\infty}$-topology on $E$ is the final topology w.r.t. all smooth curves $c:\RR\rightarrow E$. $E$ is called {\itshape convenient} if $E$ is $c^{\infty}$-complete which is equivalent for $E$ to be Mackey-complete and for $E_B$ to be a Banach space for every bounded absolutely convex subset $B$ of $E$. We refer to \cite{KM97} or the appendix in \cite{KMRc} for more details and proofs.\vspace{6pt}

{\itshape Convention:} Let $\star\in\{M,\omega,\mathcal{M}\}$, then write $\mathcal{E}_{[\star]}$ if either $\mathcal{E}_{\{\star\}}$ or $\mathcal{E}_{(\star)}$ is considered, but not mixing the cases if statements involve more than one $\mathcal{E}_{[\star]}$ symbol. The same notation will be used for the conditions, so write $(\mathcal{M}_{[\star]})$ for either $(\mathcal{M}_{\{\star\}})$ or $(\mathcal{M}_{(\star)})$.

\section{Basic definitions}
\subsection{Weight sequences and classes of ultradifferentiable functions $\mathcal{E}_{[M]}$}
A weight sequence is an arbitrary sequence of positive real numbers $M=(M_k)_k\in\RR_{>0}^{\NN}$. We introduce also $m=(m_k)_k$ defined by $m_k:=\frac{M_k}{k!}$ and $\mu_k:=\frac{M_k}{M_{k-1}}$, $\mu_0:=1$. $M$ is called {\itshape normalized} if $1=M_0\le M_1$ holds.

$(1)$ $M$ is {\itshape log-convex} if
$$\hypertarget{lc}{(\text{lc})}:\Leftrightarrow\;\forall\;j\in\NN:\;M_j^2\le M_{j-1} M_{j+1}.$$
$M$ is log-convex if and only if $(\mu_k)_k$ is increasing. If $M$ is log-convex and $M_0=1$, then
$$\hypertarget{alg}{(\text{alg})}\;:\Leftrightarrow\;\exists\;C\ge 1\;\forall\;j,k\in\NN:\;M_j M_k\le C^{j+k} M_{j+k}$$
holds with $C=1$ and the mapping $j\mapsto(M_j)^{1/j}$ is increasing, see e.g. \cite[Lemma 2.0.4, Lemma 2.0.6]{diploma}.

$M$ is called {\itshape strongly log-convex} if
$$\hypertarget{slc}{(\text{slc})}:\Leftrightarrow\;\forall\;j\in\NN:\;m_j^2\le m_{j-1} m_{j+1}.$$
This condition implies \hyperlink{lc}{$(\text{lc})$} and was a basic assumptions for $M$ in \cite{KMRc}, \cite{KMRq} and \cite{KMRplot}. It guarantees all stability properties in \cite[Theorems 5,6]{characterizationstabilitypaper} for the case $\mathcal{M}=\{M\}$, see also \cite[Theorem 3.2.]{compositionpaper}. Related to this is the weaker condition
$$\hypertarget{FdB}{(\text{FdB})}:\Leftrightarrow\;\exists\;D\ge 1\;\forall\;k\in\NN:\;\;m^{\circ}_k\le D^k m_k,$$
which is called the {\itshape Faà-di-Bruno-property}, see \cite[3.3.]{compositionpaper}. For $m^{\circ}=(m^{\circ}_k)_k$ we have put
$$m^{\circ}_k:=\max\{m_j m_{\alpha_1}\cdots m_{\alpha_j}:\;\alpha_i\in\NN_{>0}, \sum_{i=1}^j\alpha_i=k\},\hspace{20pt} m^{\circ}_0:=1.$$
Strongly log-convexity is also related to
$$\hypertarget{rai}{(\text{rai})}:\Leftrightarrow\;\exists\;C\ge 1\;\forall\;1\le j\le k:\;(m_j)^{1/j}\le C(m_k)^{1/k},$$
see \cite{compositionpaper} and \cite{characterizationstabilitypaper}.

$(2)$ $M$ has {\itshape moderate growth} if
$$\hypertarget{mg}{(\text{mg})}:\Leftrightarrow\exists\;C\ge 1\;\forall\;j,k\in\NN:\;M_{j+k}\le C^{j+k} M_j M_k.$$
This condition implies {\itshape derivation closedness}:
$$\hypertarget{dc}{(\text{dc})}\;:\Leftrightarrow\exists\;C\ge 1\;\forall\;j\in\NN:\;M_{j+1}\le C^{j+1} M_j.$$
In both conditions one can replace the sequence $M$ by $m$.

$(3)$ For $M=(M_p)_p$ and $N=(N_p)_p$ we write $M\le N$ if and only if $M_p\le N_p$ for all $p\in\NN$. Moreover we define
$$M\hypertarget{mpreceq}{\preceq}N:\Leftrightarrow\;\exists\;C_1,C_2\ge 1\;\forall\;j\in\NN:\; M_j\le C_2 C_1^j N_j\Longleftrightarrow\sup_{p\in\NN_{>0}}\left(\frac{M_p}{N_p}\right)^{1/p}<+\infty$$
and we call the sequences {\itshape equivalent} if
$$M\hypertarget{approx}{\approx}N:\Leftrightarrow\;M\hyperlink{mpreceq}{\preceq}N\;\text{and}\;N\hyperlink{mpreceq}{\preceq}M.$$
We will write
$$M\hypertarget{mtriangle}{\vartriangleleft}N:\Leftrightarrow\;\forall\;h>0\;\exists\;C_h\ge 1\;\forall\;j\in\NN:\; M_j\le C_h h^j N_j\Longleftrightarrow\lim_{p\rightarrow\infty}\left(\frac{M_p}{N_p}\right)^{1/p}=0.$$
For convenience we introduce the following set:
$$\hypertarget{LCset}{\mathcal{LC}}:=\{M\in\RR_{>0}^{\NN}:\;M\;\text{is normalized, log-convex},\;\lim_{k\rightarrow\infty}(M_k)^{1/k}=+\infty\}.$$
Let $r,s\in\NN_{>0}$ and $U\subseteq\RR^r$ be non-empty open. We introduce the ultradifferentiable class of Roumieu-type by
$$\mathcal{E}_{\{M\}}(U,\RR^s):=\{f\in\mathcal{E}(U,\RR^s):\;\forall\;K\subseteq U\;\text{compact}\;\exists\;h>0:\;\|f\|_{M,K,h}<+\infty\},$$
and the class of Beurling-type by
$$\mathcal{E}_{(M)}(U,\RR^s):=\{f\in\mathcal{E}(U,\RR^s):\;\forall\;K\subseteq U\;\text{compact}\;\forall\;h>0:\;\|f\|_{M,K,h}<+\infty\},$$
where we have put
\begin{equation}\label{semi-norm-2}
\|f\|_{M,K,h}:=\sup_{k\in\NN,x\in K}\frac{\|f^{(k)}(x)\|_{L^k(\RR^r,\RR^s)}}{h^{k} M_k}.
\end{equation}
For compact sets $K$ with smooth boundary
$$\mathcal{E}_{M,h}(K,\RR^s):=\{f\in\mathcal{E}(K,\RR^s): \|f\|_{M,K,h}<+\infty\}$$
is a Banach space and we have the topological vector space representations
\begin{equation}\label{repr3}
\mathcal{E}_{\{M\}}(U,\RR^s):=\underset{K\subseteq U}{\varprojlim}\;\underset{h>0}{\varinjlim}\;\mathcal{E}_{M,h}(K,\RR^s)=\underset{K\subseteq U}{\varprojlim}\;\mathcal{E}_{\{M\}}(K,\RR^s)
\end{equation}
and
\begin{equation}\label{repr4}
\mathcal{E}_{(M)}(U,\RR^s):=\underset{K\subseteq U}{\varprojlim}\;\underset{h>0}{\varprojlim}\;\mathcal{E}_{M,h}(K,\RR^s)=\underset{K\subseteq U}{\varprojlim}\;\mathcal{E}_{(M)}(K,\RR^s).
\end{equation}

We recall some facts for log-convex $M$:

\begin{itemize}
\item[$(i)$] Put $\mathcal{E}^{\text{global}}_{\{M\}}(U,\RR^s):=\{f\in\mathcal{E}(U,\RR^s):\;\exists\;h>0\;\|f\|_{M,U,h}<+\infty\}$. There exist {\itshape characteristic functions}
$$\hypertarget{chf}{(\text{chf})}\;:\Leftrightarrow\exists\;\theta_M\in\mathcal{E}^{\text{global}}_{\{M\}}(\RR,\RR): \forall\;j\in\NN: \left|\theta_M^{(j)}(0)\right|\ge M_j,$$
and $\tilde{\theta}_M\in\mathcal{E}^{\text{global}}_{\{M\}}(\RR,\CC)$ with
\begin{equation}\label{characcomplex}
\forall\;j\in\NN:\;\tilde{\theta}_M^{(j)}(0)=(\sqrt{-1})^j s_j,\hspace{15pt}s_j:=\sum_{k=0}^{\infty}M_k(2\mu_k)^{j-k}\ge M_j,
\end{equation}
hence $\left|\tilde{\theta}_M^{(j)}(0)\right|\ge M_j$ for all $j\in\NN$, see \cite[Lemma 2.9.]{compositionpaper} and \cite[Theorem 1]{thilliez}. Note that the Beurling-class $\mathcal{E}^{\text{global}}_{(M)}(\RR,\RR)$ cannot contain such $\theta_M$, see \cite[Proposition 3.1.2.]{diploma}.

\item[$(ii)$] If $N$ is arbitrary, then $M\hyperlink{mpreceq}{\preceq} N\Longleftrightarrow\mathcal{E}_{\{M\}}\subseteq\mathcal{E}_{\{N\}}$ and $M\hyperlink{mtriangle}{\vartriangleleft}N\Longleftrightarrow\mathcal{E}_{\{M\}}\subseteq\mathcal{E}_{(N)}$. If $M\in\hyperlink{LCset}{\mathcal{LC}}$, then $M\hyperlink{mpreceq}{\preceq} N\Longleftrightarrow\mathcal{E}_{[M]}\subseteq\mathcal{E}_{[N]}$.

\item[$(iii)$] Both classes $\mathcal{E}_{\{M\}}$ and $\mathcal{E}_{(M)}$ are closed under pointwise multiplication, see e.g. \cite[Proposition 2.0.8]{diploma}.
\end{itemize}

\subsection{Classes of ultra-differentiable functions defined by one parameter weight matrices and basic definitions}
\begin{definition}
Let $(\Lambda,\le)$ be a partially ordered set which is both up- and downward directed, $\Lambda=\RR_{>0}$ will be the most important example. A {\itshape weight matrix} $\mathcal{M}$ associated to $\Lambda$ is a family of weight sequences $\mathcal{M}:=\{M^x\in\RR_{>0}^{\NN}: x\in\Lambda\}$ such that
$$\hypertarget{Marb}{(\mathcal{M})}:\Leftrightarrow\;\forall\;x\in\Lambda:\;M^x\;\text{is normalized, increasing},\;M^x\le M^y\;\text{for}\;x\le y.$$
We call $\mathcal{M}$ {\itshape standard log-convex,} if
$$\hypertarget{Msc}{(\mathcal{M}_{\on{sc}})}:\Leftrightarrow(\mathcal{M})\;\text{and}\;\forall\;x\in\Lambda:\;M^x\in\hyperlink{LCset}{\mathcal{LC}}.$$
Also $m^x_k:=\frac{M^x_k}{k!}$ and $\mu^x_k:=\frac{M^x_k}{M^x_{k-1}}$, $\mu^x_0:=1$, will be used.
\end{definition}

We introduce ultradifferentiable classes of Roumieu- and Beurling-type defined by $\mathcal{M}$ as follows (see also \cite[4.2.]{compositionpaper}):

Let $r,s\in\NN_{>0}$, let $U\subseteq\RR^r$ be non-empty and open. For all $K\subseteq U$ compact we put
\begin{equation}\label{generalroumieu}
\mathcal{E}_{\{\mathcal{M}\}}(K,\RR^s):=\bigcup_{x\in\Lambda}\mathcal{E}_{\{M^x\}}(K,\RR^s)\hspace{20pt}\mathcal{E}_{\{\mathcal{M}\}}(U,\RR^s):=\bigcap_{K\subseteq U}\bigcup_{x\in\Lambda}\mathcal{E}_{\{M^x\}}(K,\RR^s)
\end{equation}
and
\begin{equation}\label{generalbeurling}
\mathcal{E}_{(\mathcal{M})}(K,\RR^s):=\bigcap_{x\in\Lambda}\mathcal{E}_{(M^x)}(K,\RR^s)\hspace{20pt}\mathcal{E}_{(\mathcal{M})}(U,\RR^s):=\bigcap_{x\in\Lambda}\mathcal{E}_{(M^x)}(U,\RR^s).
\end{equation}
For a compact set $K\subseteq\RR^r$ (with smooth boundary) we have
$$\mathcal{E}_{\{\mathcal{M}\}}(K,\RR^s):=\underset{x\in\Lambda}{\varinjlim}\;\underset{h>0}{\varinjlim}\;\mathcal{E}_{M^x,h}(K,\RR^s),$$
and so for $U\subseteq\RR^r$ non-empty open
\begin{equation}\label{generalroumieu1}
\mathcal{E}_{\{\mathcal{M}\}}(U,\RR^s):=\underset{K\subseteq U}{\varprojlim}\;\underset{x\in\Lambda}{\varinjlim}\;\underset{h>0}{\varinjlim}\;\mathcal{E}_{M^x,h}(K,\RR^s),
\end{equation}
and for the Beurling-case we get
\begin{equation}\label{generalbeurling1}
\mathcal{E}_{(\mathcal{M})}(U,\RR^s):=\underset{K\subseteq U}{\varprojlim}\;\underset{x\in\Lambda}{\varprojlim}\;\underset{h>0}{\varprojlim}\;\mathcal{E}_{M^x,h}(K,\RR^s).
\end{equation}

Instead of compact sets $K$ with smooth boundary one can also consider open $K\subseteq U$ with $\overline{K}$ compact in $U$, or one can work with {\itshape Whitney jets} on compact $K$.

If $\Lambda=\RR_{>0}$ we can assume that all occurring limits are countable and so $\mathcal{E}_{(\mathcal{M})}(U,\RR^s)$ is a {\itshape Fréchet space}. Moreover $\underset{x\in\Lambda}{\varinjlim}\;\underset{h>0}{\varinjlim}\;\mathcal{E}_{M^x,h}(K,\RR^s)=\underset{n\in\NN_{>0}}{\varinjlim}\;\;\mathcal{E}_{M^n,n}(K,\RR^s)$ is a {\itshape Silva space}, i.e. a countable inductive limit of Banach spaces with compact connecting mappings. For more details concerning the locally convex topology on these spaces we refer to \cite[4.2.-4.4.]{compositionpaper}.

\subsection{Conditions for a weight matrix $\mathcal{M}=\{M^x: x\in\Lambda\}$}
We are going to introduce now some conditions on $\mathcal{M}$ which will be needed frequently, see also \cite[4.1.]{compositionpaper}.

{\itshape Roumieu-type-conditions}
\begin{itemize}
\item[\hypertarget{R-dc}{$(\mathcal{M}_{\{\text{dc}\}})$}] $\forall\;x\in\Lambda\;\exists\;C>0\;\exists\;y\in\Lambda\;\forall\;j\in\NN: M^x_{j+1}\le C^{j+1} M^y_j$
\item[\hypertarget{R-mg}{$(\mathcal{M}_{\{\text{mg}\}})$}] $\forall\;x\in\Lambda\;\exists\;C>0\;\exists\;y_1,y_2\in\Lambda\;\forall\;j,k\in\NN: M^x_{j+k}\le C^{j+k} M^{y_1}_j M^{y_2}_k$
\item[\hypertarget{R-alg}{$(\mathcal{M}_{\{\text{alg}\}})$}] $\forall\;x_1,x_2\in\Lambda\;\exists\;C>0\;\exists\;y\in\Lambda\;\forall\;j,k\in\NN: M^{x_1}_j M^{x_2}_k\le C^{j+k} M^y_{j+k}$
\item[\hypertarget{R-L}{$(\mathcal{M}_{\{\text{L}\}})$}] $\forall\;C>0\;\forall\;x\in\Lambda\;\exists\;D>0\;\exists\;y\in\Lambda\;\forall\;k\in\NN: C^k M^x_k\le D M^y_k$
\item[\hypertarget{R-strict}{$(\mathcal{M}_{\{\text{strict}\}})$}] $\forall\;x\in\Lambda\;\exists\;y\in\Lambda\;:\;\;\sup_{k\in\NN_{>0}}\left(\frac{M^y_k}{M^x_k}\right)^{1/k}=+\infty$
\item[\hypertarget{R-FdB}{$(\mathcal{M}_{\{\text{FdB}\}})$}] $\forall\;x\in\Lambda\;\exists\;y\in\Lambda:(m^x)^{\circ}\hyperlink{mpreceq}{\preceq}m^y$
\item[\hypertarget{R-rai}{$(\mathcal{M}_{\{\text{rai}\}})$}] $\forall\;x\in\Lambda\;\exists\;y\in\Lambda\;\exists\;H>0:\;(m^x_q)^{1/q}\le
H(m^y_p)^{1/p},\;1\le q\le p$
\end{itemize}

{\itshape Beurling-type-conditions}
\begin{itemize}
\item[\hypertarget{B-dc}{$(\mathcal{M}_{(\text{dc})})$}] $\forall\;x\in\Lambda\;\exists\;C>0\;\exists\;y\in\Lambda\;\forall\;j\in\NN: M^y_{j+1}\le C^{j+1} M^x_j$
\item[\hypertarget{B-mg}{$(\mathcal{M}_{(\text{mg})})$}] $\forall\;x_1,x_2\in\Lambda\;\exists\;C>0\;\exists\;y\in\Lambda\;\forall\;j,k\in\NN: M^y_{j+k}\le C^{j+k} M^{x_1}_j M^{x_2}_k$
\item[\hypertarget{B-alg}{$(\mathcal{M}_{(\text{alg})})$}] $\forall\;x\in\Lambda\;\exists\;C>0\;\exists\;y_1,y_2\in\Lambda\;\forall\;j,k\in\NN: M^{y_1}_j M^{y_2}_k\le C^{j+k} M^x_{j+k}$
\item[\hypertarget{B-L}{$(\mathcal{M}_{(\text{L})})$}] $\forall\;C>0\;\forall\;x\in\Lambda\;\exists\;D>0\;\exists\;y\in\Lambda\;\forall\;k\in\NN: C^k M^y_k\le D M^x_k$
\item[\hypertarget{B-strict}{$(\mathcal{M}_{(\text{strict})})$}] $\forall\;x\in\Lambda\;\exists\;y\in\Lambda\;:\;\;\sup_{k\in\NN_{>0}}\left(\frac{M^x_k}{M^y_k}\right)^{1/k}=+\infty$
\item[\hypertarget{B-FdB}{$(\mathcal{M}_{(\text{FdB})})$}] $\forall\;x\in\Lambda\;\exists\;y\in\Lambda: (m^y)^{\circ}\hyperlink{mpreceq}{\preceq}m^x$
\item[\hypertarget{B-rai}{$(\mathcal{M}_{(\text{rai})})$}] $\forall\;x\in\Lambda\;\exists\;y\in\Lambda\;\exists\;H>0:\;(m^y_q)^{1/q}\le H(m^x_p)^{1/p},\;1\le q\le p$
\end{itemize}

\subsection{Inclusion relations of weight matrices}
Let two matrices $\mathcal{M}=\{M^x: x\in\Lambda\}$ and $\mathcal{N}=\{N^x: x\in\Lambda'\}$ be given, then we write
$$\mathcal{M}\hypertarget{Mroumpreceq}{\{\preceq\}}\mathcal{N}\;:\Leftrightarrow\forall\;x\in\Lambda\;\exists\;y\in\Lambda':\;M^x\hyperlink{mpreceq}{\preceq}N^{y}$$
$$\mathcal{M}\hypertarget{Mbeurpreceq}{(\preceq)}\mathcal{N}\;:\Leftrightarrow\forall\;y\in\Lambda'\;\exists\;x\in\Lambda:\;M^x\hyperlink{mpreceq}{\preceq}N^{y},$$
and
$$\mathcal{M}\hypertarget{Mroumapprox}{\{\approx\}}\mathcal{N}\;:\Leftrightarrow\mathcal{M}\hyperlink{Mroumpreceq}{\{\preceq\}}\mathcal{N}\;\text{and}\;\mathcal{N}\hyperlink{Mroumpreceq}{\{\preceq\}}\mathcal{M}$$
respectively
$$\mathcal{M}\hypertarget{Mbeurapprox}{(\approx)}\mathcal{N}\;:\Leftrightarrow\mathcal{M}\hyperlink{Mbeurpreceq}{(\preceq)}\mathcal{N}\;\text{and}\;\mathcal{N}\hyperlink{Mbeurpreceq}{(\preceq)}\mathcal{M}.$$
By definition $\mathcal{M}[\preceq]\mathcal{N}$ implies $\mathcal{E}_{[\mathcal{M}]}\subseteq\mathcal{E}_{[\mathcal{N}]}$. Moreover write
$$\mathcal{M}\hypertarget{Mtriangle}{\vartriangleleft}\mathcal{N}\;:\Leftrightarrow\forall\;x\in\Lambda\;\forall\;y\in\Lambda': M^x\hyperlink{mtriangle}{\vartriangleleft}N^{y},$$
so $\mathcal{M}\hypertarget{Mtriangle}{\vartriangleleft}\mathcal{N}$ implies $\mathcal{E}_{\{\mathcal{M}\}}\subseteq\mathcal{E}_{(\mathcal{N})}$. In \cite[Proposition 4.6.]{compositionpaper} the above relations are characterized for \hyperlink{Msc}{$(\mathcal{M}_{\on{sc}})$}-matrices with $\Lambda=\Lambda'=\RR_{>0}$. In this context we introduce

$\hypertarget{R-Comega}{(\mathcal{M}_{\{\mathcal{C}^{\omega}\}})}\hspace{5pt}\exists\;x\in\Lambda:\;\liminf_{k\rightarrow\infty}(m^x_k)^{1/k}>0$,

$\hypertarget{holom}{(\mathcal{M}_{\mathcal{H}})}\hspace{16pt}\forall\;x\in\Lambda:\;\liminf_{k\rightarrow\infty}(m^x_k)^{1/k}>0$,

$\hypertarget{B-Comega}{(\mathcal{M}_{(\mathcal{C}^{\omega})})}\hspace{7pt} \forall\;x\in\Lambda:\;\lim_{k\rightarrow\infty}(m^x_k)^{1/k}=+\infty$.\vspace{6pt}

If \hyperlink{R-Comega}{$(\mathcal{M}_{\{\mathcal{C}^{\omega}\}})$} holds, then $\mathcal{C}^{\omega}\subseteq\mathcal{E}_{\{\mathcal{M}\}}$, if \hyperlink{B-Comega}{$(\mathcal{M}_{(\mathcal{C}^{\omega})})$} then $\mathcal{C}^{\omega}\subseteq\mathcal{E}_{(\mathcal{M})}$. Finally if \hyperlink{holom}{$(\mathcal{M}_{\mathcal{H}})$}, then the restrictions of entire functions are contained in $\mathcal{E}_{(\mathcal{M})}$, see \cite[Proposition 4.6.]{compositionpaper}.

{\itshape Conventions:}
\begin{itemize}
\item[$(i)$] If $\Lambda=\RR_{>0}$ or $\Lambda=\NN_{>0}$, then these sets are always regarded with its natural order $\le$.
\item[$(ii)$] We will call $\mathcal{M}$ {\itshape constant} if $\mathcal{M}=\{M\}$ or more generally if $M^x\hyperlink{approx}{\approx}M^y$ for all $x,y\in\Lambda$ and which violates both \hyperlink{R-strict}{$(\mathcal{M}_{\{\on{strict}\}})$} and \hyperlink{B-strict}{$(\mathcal{M}_{(\on{strict})})$}. Otherwise it will be called {\itshape non-constant}.
\end{itemize}

\subsection{Weight functions and classes of ultradifferentiable functions $\mathcal{E}_{[\omega]}$}
A function $\omega:[0,\infty)\rightarrow[0,\infty)$ (sometimes $\omega$ is extended to $\CC$ by $\omega(x):=\omega(|x|)$) is called a {\itshape weight function} if
\begin{itemize}
\item[$\hypertarget{om0}{(\omega_0)}$] $\omega$ is continuous, on $[0,\infty)$ increasing, $\omega(x)=0$ for $x\in[0,1]$ (w.l.o.g.) and $\lim_{x\rightarrow\infty}\omega(x)=+\infty$.
\end{itemize}
Moreover we consider the following conditions:
\begin{itemize}
\item[\hypertarget{om1}{$(\omega_1)}$] $\omega(2t)=O(\omega(t))$ as $t\rightarrow+\infty$.

\item[\hypertarget{om2}{$(\omega_2)$}] $\omega(t)=O(t)$ as $t\rightarrow\infty$.

\item[\hypertarget{om3}{$(\omega_3)$}] $\log(t)=o(\omega(t))$ as $t\rightarrow+\infty$ ($\Leftrightarrow\lim_{t\rightarrow+\infty}\frac{t}{\varphi_{\omega}(t)}=0$).

\item[\hypertarget{om4}{$(\omega_4)$}] $\varphi_{\omega}:t\mapsto\omega(e^t)$ is a convex function on $\RR$.

\item[\hypertarget{om5}{$(\omega_5)$}] $\omega(t)=o(t)$ as $t\rightarrow+\infty$.

\item[\hypertarget{om6}{$(\omega_6)$}] $\exists\;H\ge 1\;\forall\;t\ge 0:\;2\omega(t)\le\omega(Ht)+H$.

\item[\hypertarget{om1pp}{$(\omega_{1'})$}] $\exists\;D>0: \exists\;t_0>0: \forall\:\lambda\ge 1: \forall\;t\ge t_0: \omega(\lambda t)\le D\lambda\omega(t)$.
\end{itemize}
An interesting example is $\omega_s(t):=\max\{0,\log(t)^s\}$, $s>1$, which satisfies all listed properties except \hyperlink{om6}{$(\omega_6)$}. For convenience we define the sets
$$\hypertarget{omset0}{\mathcal{W}_0}:=\{\omega:[0,\infty)\rightarrow[0,\infty): \omega\;\text{has}\;\hyperlink{om0}{(\omega_0)},\hyperlink{om3}{(\omega_3)},\hyperlink{om4}{(\omega_4)}\},$$
$$\hypertarget{omset1}{\mathcal{W}}:=\{\omega\in\mathcal{W}_0: \omega\;\text{has}\;\hyperlink{om1}{(\omega_1)}\}.$$
For $\omega\in\hyperlink{omset0}{\mathcal{W}_0}$ we define the {\itshape Legendre-Fenchel-Young-conjugate} $\varphi^{*}_{\omega}$ by
$$\varphi^{*}_{\omega}(x):=\sup\{xy-\varphi_{\omega}(y): y\ge 0\},\hspace{15pt}x\ge 0.$$
It is a convex increasing function, $\varphi^{*}_{\omega}(0)=0$, $\varphi^{**}_{\omega}=\varphi_{\omega}$, $\lim_{x\rightarrow\infty}\frac{x}{\varphi^{*}_{\omega}(x)}=0$ and finally $x\mapsto\frac{\varphi_{\omega}(x)}{x}$ and $x\mapsto\frac{\varphi^{*}_{\omega}(x)}{x}$ are increasing on $[0,+\infty)$, see e.g. \cite[Remark 1.3., Lemma 1.5.]{BraunMeiseTaylor90}.

For $\sigma,\tau\in\hyperlink{omset1}{\mathcal{W}}$ we write
$$\sigma\hypertarget{ompreceq}{\preceq}\tau:\Leftrightarrow\tau(t)=O(\sigma(t)),\;\;\text{as}\;t\rightarrow+\infty$$
and call them equivalent if
$$\sigma\hypertarget{sim}{\sim}\tau:\Leftrightarrow\sigma\hyperlink{ompreceq}{\preceq}\tau\;\text{and}\;\tau\hyperlink{ompreceq}{\preceq}\sigma.$$
Let $r,s\in\NN_{>0}$, $U\subseteq\RR^r$ be a non-empty open set and $\omega\in\hyperlink{omset0}{\mathcal{W}_0}$. The Roumieu-type space is defined by
$$\mathcal{E}_{\{\omega\}}(U,\RR^s):=\{f\in\mathcal{E}(U,\RR^s):\;\forall\;K\subseteq U\;\text{compact}\;\exists\;l>0:\;\|f\|_{\omega,K,l}<+\infty\}$$
and the Beurling-type space by
$$\mathcal{E}_{(\omega)}(U,\RR^s):=\{f\in\mathcal{E}(U,\RR^s):\;\forall\;K\subseteq U\;\text{compact}\;\forall\;l>0:\;\|f\|_{\omega,K,l}<+\infty\},$$
where we have put
\begin{equation}\label{semi-norm-1}
\|f\|_{\omega,K,l}:=\sup_{k\in\NN,x\in K}\frac{\|f^{(k)}(x)\|_{L^k(\RR^r,\RR^s)}}{\exp(\frac{1}{l}\varphi^{*}_{\omega}(lk))}
\end{equation}
and $f^{(k)}(x)$ denotes the $k$-th order Fréchet derivative at $x$. For compact sets $K$ with smooth boundary
$$\mathcal{E}_{\omega,l}(K,\RR^s):=\{f\in\mathcal{E}(K,\RR^s): \|f\|_{\omega,K,l}<+\infty\}$$
is a Banach space and we have the topological vector space representations
\begin{equation}\label{repr1}
\mathcal{E}_{\{\omega\}}(U,\RR^s):=\underset{K\subseteq U}{\varprojlim}\;\underset{l>0}{\varinjlim}\;\mathcal{E}_{\omega,l}(K,\RR^s)=\underset{K\subseteq U}{\varprojlim}\;\mathcal{E}_{\{\omega\}}(K,\RR^s)
\end{equation}
and
\begin{equation}\label{repr2}
\mathcal{E}_{(\omega)}(U,\RR^s):=\underset{K\subseteq U}{\varprojlim}\;\underset{l>0}{\varprojlim}\;\mathcal{E}_{\omega,l}(K,\RR^s)=\underset{K\subseteq U}{\varprojlim}\;\mathcal{E}_{(\omega)}(K,\RR^s).
\end{equation}

A new idea introduced in \cite[Chapter 5]{compositionpaper} was the following:
\begin{itemize}
\item[$(i)$] To each $\omega\in\hyperlink{omset1}{\mathcal{W}}$ we can associate a \hyperlink{Msc}{$(\mathcal{M}_{\text{sc}})$} weight matrix $\Omega=\{\Omega^l=(\Omega^l_j)_{j\in\NN}: l >0\}$ by

    \centerline{\fbox{$\Omega^l_j:=\exp\left(\frac{1}{l}\varphi^{*}_{\omega}(lj)\right)$.}}

\item[$(ii)$] $\Omega$ has always \hyperlink{R-mg}{$(\mathcal{M}_{\{\text{mg}\}})$} and \hyperlink{B-mg}{$(\mathcal{M}_{(\text{mg})})$}, \hyperlink{R-L}{$(\mathcal{M}_{\{\text{L}\}})$} and \hyperlink{B-L}{$(\mathcal{M}_{(\text{L})})$}. If $\omega$ is sub-additive, then  \hyperlink{R-FdB}{$(\mathcal{M}_{\{\text{FdB}\}})$} and \hyperlink{B-FdB}{$(\mathcal{M}_{(\text{FdB})})$} hold, see \cite[Lemma 6.1.]{compositionpaper}. Equivalent weight functions $\omega$ yield equivalent weight matrices w.r.t. both \hyperlink{Mbeurapprox}{$(\approx)$} and \hyperlink{Mroumapprox}{$\{\approx\}$}.
\item[$(iii)$] $\mathcal{E}_{[\Omega]}=\mathcal{E}_{[\omega]}$ holds as locally convex vector spaces, so defining classes of ultradifferentiable functions by weight matrices as in \eqref{generalroumieu} and \eqref{generalbeurling} is a common generalization of defining them by using a single weight sequence $M$, i.e. a constant weight matrix, or a weight function $\omega\in\hyperlink{omset1}{\mathcal{W}}$. But one is also able to describe classes which cannot be described neither by a weight function nor by a weight sequence, e.g. the class defined by the Gevrey-matrix $\mathcal{G}:=\{(p!^{s+1})_{p\in\NN}: s>0\}$, see \cite[5.19.]{compositionpaper}.
\end{itemize}

\section{Basic definitions for the convenient setting}\label{section2}
\subsection{Whitney jets on Banach spaces}\label{subsection21}
We recall the notation of \cite[Chapter 3]{KMRplot}. Let $E,F$ be Banach spaces, $K\subseteq E$ compact and $U\subseteq E$ open. Let $f\in\mathcal{E}(U,F)$, then we introduce the jet mapping $j^{\infty}:\mathcal{E}(U,F)\rightarrow J^{\infty}(U,F):=\prod_{k\in\NN}\mathcal{C}(U,L^k_{\sym}(E,F))$ defined by $f\mapsto j^{\infty}(f)=(f^{(k)})_{k\in\NN}$. For an arbitrary subset $X\subseteq E$ and an infinite jet $f=(f^k)_{k\in\NN}$ we introduce the {\itshape Taylor polynomial} $(T^n_y f)^k: X\rightarrow L^k_{\sym}(E,F)$ of order $n$ at the point $y$ as follows:
$$(T^n_y f)^k(x)(v_1,\dots,v_k):=\sum_{j=0}^n\frac{1}{j!}f^{j+k}(y)(x-y,\dots,x-y,v_1,\dots,v_k).$$
The {\itshape remainder} is given by
$$(R^n_yf)^k(x):=f^k(x)-(T^n_y f)^k(x)=(T^n_x f)^k(x)-(T^n_y f)^k(x)$$
and so $(R^n_yf)^k(x)\in L^k_{\sym}(E,F)$. We put now
$$\|f\|_k:=\sup\left\{\|f^k(x)\|_{L^k_{\sym}(E,F)}: x\in K\right\}$$
and
$$|\|f|\|_{n,k}:=\sup\left\{(n+1)!\frac{\|(R_y^n f)^k(x)\|_{L^k_{\sym}(E,F)}}{\|x-y\|^{n+1}}: x,y\in K, x\neq y\right\}.$$
We supply $\mathcal{E}(U,F)$ with the seminorms $f\mapsto\|j^{\infty}(f)|_K\|_k$, where $K\subseteq U$ is a compact set and $k\in\NN$. If $K\subseteq E$ is compact and {\itshape convex}, then we introduce the space $\mathcal{E}(E\supseteq K,F)$ of {\itshape Whitney-jets} on $K$ by
$$\mathcal{E}(E\supseteq K,F):=\left\{f=(f^k)_{k\in\NN}\in\prod_{k\in\NN}\mathcal{C}(K,L^k_{\sym}(E,F)): |\|f|\|_{n,k}<+\infty\;\forall\;n,k\in\NN\right\}$$
and we supply these spaces with both seminorms $\|f\|_k$ and $|\|f|\|_{n,k}$ for $k,n\in\NN$. Finally recall \cite[Lemma 3.1.]{KMRplot}:

\begin{lemma}\label{KMRfirstlemma}
Let $E$ and $F$ be Banach spaces and $K\subseteq E$ be a compact convex subset. Then $\mathcal{E}(E\supseteq K,F)$ is a Fréchet space.
\end{lemma}

\subsection{Classes of ultra-differentiable mappings defined by a weight matrix}\label{subsection22}
Let $\mathcal{M}:=\{M^x: x\in\Lambda\}$ be \hyperlink{Marb}{$(\mathcal{M})$}, $E$ and $F$ be Banach spaces and $K\subseteq E$ a compact subset. Then, as in \cite[4.1.]{KMRplot}, for $x\in\Lambda$ and $h>0$ we define
$$\mathcal{E}_{M^x,h}(E\supseteq K,F):=\left\{(f^j)_j\in\prod_{j\in\NN}\mathcal{C}(K,L^j_{\sym}(E,F)): \|f\|^J_{M^x,h}<+\infty\right\},$$
where
$$\|f\|^J_{M^x,h}:=\max\left\{\sup\left\{\frac{\|f\|_k}{h^k M^x_k}: k\in\NN\right\}, \sup\left\{\frac{|\|f|\|_{n,k}}{h^{n+k-1} M^x_{n+k+1}}: k,n\in\NN\right\}\right\}.$$
For open $U\subseteq E$ and compact $K\subseteq U$ we introduce the space
$$\mathcal{E}_{M^x,K,h}(U,F):=\left\{f\in\mathcal{E}(U,F): j^{\infty}(f)\big|_K\in\mathcal{E}_{M^x,h}(E\supseteq K,F)\right\},$$
with semi-norm $f\mapsto\left\|j^{\infty}(f)\big|_K\right\|^J_{M^x,h}$. It is not Hausdorff and for infinite dimensional $E$ its Hausdorff quotient will not always be complete. Note that if $K$ is assumed to be {\itshape convex}, then we can take on $\mathcal{E}_{M^x,K,h}(U,F)$ also the semi-norm
$$f\mapsto\sup\left\{\frac{\|f^{(n)}(a)\|_{L^n_{\sym}(E,F)}}{h^n M^x_n}: a\in K, n\in\NN\right\}=:\|f\|^J_{M^x,K,h}.$$
Thus we see that $\mathcal{E}_{M^x,K,h}(U,F)=\left\{f\in\mathcal{E}(U,F): (\|j^{\infty}(f)\big|_K\|_k)_k\in\mathcal{F}_{M^x,h}\right\}$ holds with
$$\mathcal{F}_{M^x,h}:=\left\{(f_k)_k\in\RR_{>0}^{\NN}: \exists\;C>0:\forall\;k\in\NN: |f_k|\le C h^k M^x_k\right\}.$$
The bounded sets $\mathcal{B}$ in $\mathcal{E}_{M^x,K,h}(U,F)$ are exactly those $\mathcal{B}\subseteq\mathcal{E}(U,F)$ such that $(b_m)_m\in\mathcal{F}_{M^x,h}$ with $b_m:=\sup\left\{\|j^{\infty}(f)\big|_K\|_m: f\in\mathcal{B}\right\}$.\vspace{6pt}

Let $U\subseteq E$ be convex open and $K\subseteq U$ be convex compact, then define
$$\mathcal{E}_{(\mathcal{M})}(E\supseteq K,F):=\underset{x\in\Lambda,h>0}{\varprojlim}\;\mathcal{E}_{M^x,h}(E\supseteq K,F)$$
$$\mathcal{E}_{\{\mathcal{M}\}}(E\supseteq K,F):=\underset{x\in\Lambda,h>0}{\varinjlim}\;\mathcal{E}_{M^x,h}(E\supseteq K,F)$$
and finally
\begin{equation}\label{plotfirstdef}
\mathcal{E}_{[\mathcal{M}]}(U,F):=\underset{K\subseteq U}{\varprojlim}\;\mathcal{E}_{[\mathcal{M}]}(E\supseteq K,F),
\end{equation}
i.e.
$$\mathcal{E}_{[\mathcal{M}]}(U,F):=\left\{f\in\mathcal{E}(U,F): \forall\;K: (f^{(k)}\big|_K)\in\mathcal{E}_{[\mathcal{M}]}(E\supseteq K,F)\right\},$$
where $K$ runs through all compact and convex subsets of $U$.

If $\Lambda=\RR_{>0}$, then we can restrict in both cases to the countable diagonal, see also \cite[4.2.-4.4.]{compositionpaper}. We have $\mathcal{E}_{(\mathcal{M})}(E\supseteq K,F)=\underset{n\in\NN_{>0}}{\varprojlim}\;\mathcal{E}_{M^{1/n},1/n}(E\supseteq K,F)$ and $\mathcal{E}_{\{\mathcal{M}\}}(E\supseteq K,F)=\underset{n\in\NN_{>0}}{\varinjlim}\;\mathcal{E}_{M^n,n}(E\supseteq K,F)$.

As already mentioned in \cite[Proposition 4.1. $(3)$]{KMRplot} the space $\mathcal{E}_{\{\mathcal{M}\}}(E\supseteq K,F)$ is not a {\itshape Silva space} for infinite dimensional $E$, because the connecting mappings in the inductive limit $\underset{x\in\Lambda,h>0}{\varinjlim}\;\mathcal{E}_{M^x,h}(E\supseteq K,F)$ are not compact any more. The set $\mathcal{B}:=\{\alpha\in E^{'}: \|\alpha\|\le 1\}$ is bounded in $\mathcal{E}_{M^k,k}(E\supseteq K,\RR)$ for each $k\ge 1$. We have $\|\alpha\|_0=\sup\{|\alpha(x)|: x\in K\}\le\sup\{\|x\|: x\in K\}$, $\|\alpha\|_1=\|\alpha\|\le 1$ and $\|\alpha\|_m=0$ for each $m\ge 2$. Moreover $(R^n_y\alpha)^k=0$ for $n+k\ge 1$ and $(R^0_y\alpha)^0=\alpha(x-y)$. But $\mathcal{B}$ is not relatively compact in any $\mathcal{E}_{M^k,k}(E\supseteq K,\RR)$, $k\ge 1$, because it is not even pointwise relatively compact in $\mathcal{C}(K,L(E,\RR))$.\vspace{6pt}

Moreover we define
$$\mathcal{E}_{(\mathcal{M}),K}(U,F):=\underset{x\in\Lambda,h>0}{\varprojlim}\;\mathcal{E}_{M^x,K,h}(U,F)$$
$$\mathcal{E}_{\{\mathcal{M}\},K}(U,F):=\underset{x\in\Lambda,h>0}{\varinjlim}\;\mathcal{E}_{M^x,K,h}(U,F)$$
and so
$$\mathcal{E}_{(\mathcal{M}),K}(U,F)=\left\{f\in\mathcal{E}(U,F): (\|j^{\infty}(f)\big|_K\|_k)_k\in\mathcal{F}_{(\mathcal{M})}\right\}$$
$$\mathcal{E}_{\{\mathcal{M}\},K}(U,F)=\left\{f\in\mathcal{E}(U,F): (\|j^{\infty}(f)\big|_K\|_k)_k\in\mathcal{F}_{\{\mathcal{M}\}}\right\}$$
with $\mathcal{F}_{(\mathcal{M})}=\bigcap_{x\in\Lambda,h>0}\mathcal{F}_{M^x,h}$, $\mathcal{F}_{\{\mathcal{M}\}}=\bigcup_{x\in\Lambda,h>0}\mathcal{F}_{M^x,h}$.

The bounded sets $\mathcal{B}\subseteq\mathcal{E}_{[\mathcal{M}],K}(U,F)$ are exactly those $\mathcal{B}\subseteq\mathcal{E}(U,F)$ for which the sequence $(b_m)_m$, $b_m:=\sup\left\{\|j^{\infty}(f)\big|_K\|_m: f\in\mathcal{B}\right\}$, belongs to $\mathcal{F}_{[\mathcal{M}]}$.

Finally we introduce
$$\underset{K\subseteq U}{\varprojlim}\;\mathcal{E}_{[\mathcal{M}],K}(U,F)=\left\{f\in\mathcal{E}(U,F): \forall\;K: (\|j^{\infty}(f)\big|_K\|_m)_m\in\mathcal{F}_{[\mathcal{M}]}\right\}.$$

The next result generalizes \cite[Proposition 4.1.]{KMRplot}.

\begin{proposition}\label{matrixcompletness}
Let $\mathcal{M}$ be \hyperlink{Marb}{$(\mathcal{M})$} with $\Lambda=\RR_{>0}$, then the following completeness properties are valid:

\begin{itemize}
\item[$(1)$] $\mathcal{E}_{M^x,h}(E\supseteq K,F)$ is a Banach space.

\item[$(2)$] $\mathcal{E}_{(\mathcal{M})}(E\supseteq K,F)$ is a Fréchet space.

\item[$(3)$] $\mathcal{E}_{\{\mathcal{M}\}}(E\supseteq K,F)$ is a compactly regular $(LB)$-space, i.e. compact subsets are contained and compact in some step and so ($c^{\infty}$)-complete, webbed and ultrabornological.

\item[$(4)$] $\mathcal{E}_{(\mathcal{M})}(U,F)$ and $\mathcal{E}_{\{\mathcal{M}\}}(U,F)$ are complete.

\item[$(5)$] As locally convex vector spaces we have
$$\mathcal{E}_{(\mathcal{M})}(U,F)=\underset{K\subseteq U}{\varprojlim}\;\mathcal{E}_{(\mathcal{M})}(E\supseteq K,F)=\underset{K\subseteq U}{\varprojlim}\;\mathcal{E}_{(\mathcal{M}),K}(U,F)$$
and
$$\mathcal{E}_{\{\mathcal{M}\}}(U,F)=\underset{K\subseteq U}{\varprojlim}\;\mathcal{E}_{\{\mathcal{M}\}}(E\supseteq K,F)=\underset{K\subseteq U}{\varprojlim}\;\mathcal{E}_{\{\mathcal{M}\},K}(U,F).$$
\end{itemize}
\end{proposition}

\demo{Proof}
$(1)$ This was already shown in \cite[Proposition 4.1. $(1)$]{KMRplot}.

$(2)$ Holds since $\Lambda=\RR_{>0}$.

$(3)$ We can restrict to $\Lambda=\NN_{>0}$ and proceed analogously as in \cite[Proposition 4.1. $(3)$]{KMRplot}. To show that the inductive limit is compactly regular it suffices to show that there exists a sequence of increasing $0$-neighborhoods $U_n\in\mathcal{E}_{M^n,n}(E\supseteq K,F)$ such that for each $n\in\NN$ there exists $l\in\NN$ with $l\ge n$ and for which the topologies of $\mathcal{E}_{M^l,l}(E\supseteq K,F)$ and of $\mathcal{E}_{M^k,k}(E\supseteq K,F)$ coincide on $U_n$ for all $k\ge l$.\vspace{6pt}

In general, for indices $x_1\ge x_2$ and positive real numbers $h_1\ge h_2$ we have clearly by definition $\|\cdot\|^J_{M^{x_1},h_1}\le\|\cdot\|^J_{M^{x_2},h_2}$. Consider now the $\varepsilon$-Ball $U^{x,h}_{\varepsilon}(f):=\{g: \|g-f\|^J_{M^x,h}\le\varepsilon\}$ in $\mathcal{E}_{M^x,h}(E\supseteq K,F)$ and we restrict to the diagonal $x=h=n$ and identify $U^{n,n}$ with $U^n$.

We show that for arbitrary $n\in\NN_{>0}$ and $n_2>n_1:=2n$, for each $\varepsilon>0$ and $f\in U_1^{n}(0)$ there exists $\delta>0$ such that $U_{\delta}^{n_2}(f)\cap U_1^{n}(0)\subseteq U_{\varepsilon}^{n_1}(f)$.

By assumption $f\in U_1^n(0)=U_1^{n,n}(0)$ we have $\|f\|_a\le n^a M^n_a$ and $|\|f\||_{a,b}\le n^{a+b+1} M^n_{a+b+1}$ for all $a,b\in\NN$. Consider $g\in U_{\delta}^{n_2}(f)\cap U_1^{n}(0)=U_{\delta}^{n_2,n_2}(f)\cap U_1^{n,n}(0)$, then $\|g\|_a\le n^a M^n_a$, $|\|g\||_{a,b}\le n^{a+b+1} M^n_{a+b+1}$ and moreover $\|g-f\|_a\le\delta n_2^a M^{n_2}_a$, $|\|g-f\||_{a,b}\le\delta n_2^{a+b+1} M^{n_2}_{a+b+1}$ for all $a,b\in\NN$. We estimate similarly as in \cite[Proposition 4.1. $(3)$]{KMRplot}. So for given $\varepsilon>0$ consider $N\in\NN$ (minimal) with $\frac{1}{2^N}<\frac{\varepsilon}{2}$ and put $\delta:=\varepsilon\left(\frac{n_1}{n_2}\right)^{N-1}\frac{1}{M_N^{n_2}}$.\vspace{6pt}

For $a\ge N$ we have $\frac{1}{2^a}\le\frac{1}{2^N}<\frac{\varepsilon}{2}$ $(\star)$, so use triangle-inequality to get
$$\|g-f\|_a\le \|g\|_a+\|f\|_a\le 2 n^a M^n_a=2 n_1^a M^n_a\frac{1}{2^a}\underbrace{<}_{(\star)}\varepsilon n_1^a M^n_a\le\varepsilon n_1^a M^{n_1}_a$$
and the last inequality holds since $n_1=2n>n$ and so $M^n_a\le M^{n_1}_a$ for all $a\in\NN$. For $a<N$ we have
$$\|g-f\|_a\le\delta n_2^a M_a^{n_2}\le\varepsilon n_1^a\frac{M_a^{n_2}}{M_N^{n_2}}\le\varepsilon n_1^a\le\varepsilon n_1^a M^{n_1}_a,$$
because $M_a^n\le M_N^n$, $\left(\frac{n_1}{n_2}\right)^{N-1}\le\left(\frac{n_1}{n_2}\right)^{a}$ since $a<N$, $\frac{n_1}{n_2}<1$ and finally $M^{n_1}_a\ge 1$.

Analogously we can use the same estimates for $|\|\cdot\||_{a,b}$ instead of $\|\cdot\|_a$ for each $a,b\in\NN$.\vspace{6pt}

$(4)$ In the Beurling-case we have a projective limit of Fréchet spaces, in the Roumieu-case a projective limit of $(LB)$-spaces, which are all compactly regular by $(3)$ and so complete, too. Since projective limits of complete spaces are complete we are done.

$(5)$ This holds precisely by the same proof as given in \cite[Proposition 4.1. $(5)$]{KMRplot}
\qed\enddemo

Let $E,F$ be convenient, $U\subseteq E$ be $c^{\infty}$-open, then define
\begin{align*}
\mathcal{E}^{\on{b}}_{(\mathcal{M})}(U,F):=&
\Big\{f\in\mathcal{E}(U,F): \forall\;B:\forall\;K\subseteq U\cap E_B:\forall\;x\in\Lambda\;\forall\;h>0:
\\&
\Big\{\frac{f^{(k)}(a)(v_1,\dots,v_k)}{h^k M^x_k}: k\in\NN, a\in K, \|v_i\|_B\le 1\Big\}\text{is bounded in}\;F\Big\}
\\&
=\Big\{f\in\mathcal{E}(U,F): \forall\;B:\forall\;K\subseteq U\cap E_B:\forall\;x\in\Lambda\;\forall\;h>0:
\\&
\Big\{\frac{d^k_v f(a)(v_1,\dots,v_k)}{h^k M^x_k}: k\in\NN, a\in K, \|v_i\|_B\le 1\Big\}\text{is bounded in}\;F\Big\}.
\end{align*}
and
\begin{align*}
\mathcal{E}^{\on{b}}_{\{\mathcal{M}\}}(U,F):=&
\Big\{f\in\mathcal{E}(U,F): \forall\;B:\forall\;K\subseteq U\cap E_B:\exists\;x\in\Lambda\;\exists\;h>0:
\\&
\Big\{\frac{f^{(k)}(a)(v_1,\dots,v_k)}{h^k M^x_k}: k\in\NN, a\in K, \|v_i\|_B\le 1\Big\}\text{is bounded in}\;F\Big\}
\\&
=\Big\{f\in\mathcal{E}(U,F): \forall\;B:\forall\;K\subseteq U\cap E_B:\exists\;x\in\Lambda\;\exists\;h>0:
\\&
\Big\{\frac{d^k_v f(a)(v_1,\dots,v_k)}{h^k M^x_k}: k\in\NN, a\in K, \|v_i\|_B\le 1\Big\}\text{is bounded in}\;F\Big\}.
\end{align*}

$B$ runs through all closed absolutely convex bounded subsets in $E$, $E_B$ is the complete vector space generated by $B$ with the {\itshape Minkowski-functional} $\|\cdot\|_B$. Finally $K$ runs through all sets in $U\cap E_B$ which are compact w.r.t. the norm $\|\cdot\|_B$. If $E$ and $F$ both are Banach spaces and $U\subseteq E$ open we have $\mathcal{E}^{\on{b}}_{[\mathcal{M}]}(U,F)=\mathcal{E}_{[\mathcal{M}]}(U,F)$, where the latter space is introduced in \eqref{plotfirstdef}.

Now we give the most important definition:

\centerline{\fbox{$\mathcal{E}_{[\mathcal{M}]}(U,F):=\left\{f\in\mathcal{E}(U,F): \forall\;\alpha\in F^{*}:\forall\;B: \alpha\circ f\circ i_B\in\mathcal{E}_{[\mathcal{M}]}(U_B,\RR)\right\},$}}

where $B$ is running again through all closed absolutely convex bounded subsets in $E$, the mapping $i_B: E_B\rightarrow E$ denotes the inclusion of $E_B$ in $E$ and we write $U_B:=i_B^{-1}(U)$. The initial locally convex structure is now induced by all linear mappings

\centerline{\fbox{$\mathcal{E}_{[\mathcal{M}]}(i_B,\alpha): \mathcal{E}_{[\mathcal{M}]}(U,F)\longrightarrow\mathcal{E}_{[\mathcal{M}]}(U_B,\RR), \hspace{20pt} f\mapsto\alpha\circ f\circ i_B.$}}

$\mathcal{E}_{[\mathcal{M}]}(U,F)\subseteq\prod_{\alpha,B}\mathcal{E}_{[\mathcal{M}]}(U_B,\RR)$ are convenient vector spaces as $c^{\infty}$-closed subspaces in the product: Smoothness can be tested by composing with inclusions $E_B\rightarrow E$ and $\alpha\in F^{\ast}$ as mentioned in \cite[2.14.4, 1.8]{KM97}. Hence we obtain the representation
\begin{equation}\label{Banach-plots-importantdefinition2}
\mathcal{E}_{[\mathcal{M}]}(U,F):=\{f\in F^U:\;\forall\;\alpha\in F^{\ast}\;\forall\;B: \alpha\circ f\circ i_B\in\mathcal{E}_{[\mathcal{M}]}(U_B,\RR)\}.
\end{equation}
All definitions given here are clearly generalizations of the definitions in \cite[4.2.]{KMRplot} for constant matrices.

\section{Projective descriptions for $\mathcal{E}_{[\mathcal{M}]}$}\label{section3}
In this section we are going to study one of the most important new techniques in this work. Using abstract families of sequences of positive real numbers we prove projective representations for the Roumieu-class $\mathcal{E}_{\{\mathcal{M}\}}$. This technique is very important since we want to get rid of {\itshape both} existence quantifiers in the definitions of $\mathcal{E}_{\{\mathcal{M}\}}$ so we want to generalize \cite[Lemma 4.6.]{KMRplot}. Furthermore we are going to prove analogous results for the Beurling-case $\mathcal{E}_{(\mathcal{M})}$ and generalize \cite[Lemma 4.5.]{KMRplot}. To do so we have to show variations and generalizations of \cite[Lemma 9.2.]{KM97} (for the Roumieu-case) and of the Lemma between Lemma $4.5.$ and Lemma $4.6.$ in \cite{KMRplot} (for the Beurling-case).

We will obtain different projective representations for $\mathcal{E}_{[\mathcal{M}]}$. The choice of the appropriate representation depends on the application in the proofs. To show closedness under composition in section \ref{section4}, see Theorem \ref{Banach-plot-theorem} and Theorem \ref{Banach-plots-category}, we will have to use the versions using the Faà-di-Bruno-property $(\mathcal{M}_{[\on{FdB}]})$. For the exponential laws in section \ref{section5} the versions only assuming \hyperlink{Marb}{$(\mathcal{M})$} or \hyperlink{Msc}{$(\mathcal{M}_{\on{sc}})$} for $\mathcal{M}$ are sufficient.

First we have to introduce several classes of sequences of positive real numbers $(r_k)_k$ and $(s_k)_k$. It is no restriction to assume $r_0=1$ resp. $s_0=1$ (normalization) for all occurring sequences.\vspace{6pt}

$\hypertarget{r-roum}{\mathcal{R}_{\text{Roum}}}\hspace{14.5pt}:=\{(r_k)_k\in\RR_{>0}^{\NN}: r_k t^k\rightarrow 0\;\text{as}\;k\rightarrow\infty\;\;\text{for each}\;t>0\}$

$\hypertarget{r-roum-sub}{\mathcal{R}_{\text{Roum},\text{sub}}}:=\{(r_k)_k\in\mathcal{R}_{\text{Roum}}: r_{j+k}\le r_k r_j\;\forall\;j,k\in\NN\}$

$\hypertarget{r-beur}{\mathcal{R}_{\text{Beur}}}\hspace{19pt}:=\{(r_k)_k\in\RR_{>0}^{\NN}: r_k t^k\rightarrow 0\;\text{as}\;k\rightarrow\infty\;\;\text{for some}\;t>0\}$

$\hypertarget{r-beur-sub}{\mathcal{R}_{\text{Beur},\text{sub}}}\hspace{4.5pt}:=\{(r_k)_k\in\mathcal{R}_{\text{Beur}}: r_{j+k}\le r_k r_j\;\forall\;j,k\in\NN\}$

$\hypertarget{s-roum}{\mathcal{S}^{\mathcal{M}}_{\text{Roum}}}\hspace{17pt}:=\{(s_k)_k\in\RR_{>0}^{\NN}: \forall\;x\in\Lambda\;\exists\;C_x>0\;\forall\;k\in\NN: s_k m^x_k\le C_x^k\}$

$\hypertarget{S-roum}{\tilde{\mathcal{S}}^{\mathcal{M}}_{\text{Roum}}}\hspace{17pt}:=\{(s_k)_k\in\RR_{>0}^{\NN}: \forall\;x\in\Lambda\;\exists\;C_x>0\;\forall\;k\in\NN: s_k M^x_k\le C_x^k\}$

$\hypertarget{S-roum-sub}{\tilde{\mathcal{S}}^{\mathcal{M}}_{\text{Roum},\text{sub}}}\hspace{2.5pt}:=\{(s_k)_k\in\tilde{\mathcal{S}}^{\mathcal{M}}_{\text{Roum}}:\exists\;D>0\;\forall\;j,k\in\NN: s_{j+k}\le D^{j+k} s_j s_k\}$

$\hypertarget{s-roum-FdB}{\mathcal{S}^{\mathcal{M}}_{\text{Roum},\text{FdB}}}\hspace{0pt}:=\{(s_k)_k\in\mathcal{S}^{\mathcal{M}}_{\text{Roum}}:\;\exists\;(\hat{s}_k)_k\in\mathcal{S}^{\mathcal{M}}_{\text{Roum}}\;\exists\;D>0\;\forall\;k\in\NN: s_k\le D^k(\hat{s}_o)_k\}$

$\hypertarget{s-beur}{\mathcal{S}^{\mathcal{M}}_{\text{Beur}}}\hspace{21pt}:=\{(s_k)_k\in\RR_{>0}^{\NN}: \exists\;x\in\Lambda\;\exists\;C_x>0\;\forall\;k\in\NN: s_k m^x_k\le C_x^k\}$

$\hypertarget{S-beur}{\tilde{\mathcal{S}}^{\mathcal{M}}_{\text{Beur}}}\hspace{21pt}:=\{(s_k)_k\in\RR_{>0}^{\NN}: \exists\;x\in\Lambda\;\exists\;C_x>0\;\forall\;k\in\NN: s_k M^x_k\le C_x^k\}$

$\hypertarget{S-beur-sub}{\tilde{\mathcal{S}}^{\mathcal{M}}_{\text{Beur},\text{sub}}}\hspace{6.5pt}:=\{(s_k)_k\in\tilde{\mathcal{S}}^{\mathcal{M}}_{\text{Beur}}:\exists\;D>0\;\forall\;j,k\in\NN: s_{j+k}\le D^{j+k} s_j s_k\}$

$\hypertarget{s-beur-FdB}{\mathcal{S}^{\mathcal{M}}_{\text{Beur},\text{FdB}}}\hspace{4pt}:=\{(s_k)_k\in\mathcal{S}^{\mathcal{M}}_{\text{Beur}}:\;\exists\;(\hat{s}_k)_k\in\mathcal{S}^{\mathcal{M}}_{\text{Beur}}\;\exists\;D>0\;\forall\;k\in\NN: s_k\le D^k(\hat{s}_o)_k\}$
\vspace{6pt}

For $(s_k)_k\in\mathcal{S}^{\mathcal{M}}_{\text{Roum}},\mathcal{S}^{\mathcal{M}}_{\text{Beur}}$ we have put
$$(s_o)_k:=\min\{s_j s_{\alpha_1}\cdots s_{\alpha_j}: \alpha_i\in\NN_{>0}, \alpha_1+\dots+\alpha_j=k\},\;\;\;(s_o)_0:=1.$$
By definition $\hyperlink{S-roum}{\tilde{\mathcal{S}}^{\mathcal{M}}_{\text{Roum}}}\subseteq\hyperlink{s-roum}{\mathcal{S}^{\mathcal{M}}_{\text{Roum}}}$ and $(s_k)_k\in\hyperlink{S-roum}{\tilde{\mathcal{S}}^{\mathcal{M}}_{\text{Roum}}}$ if and only if $(k!s_k)_k\in\hyperlink{s-roum}{\mathcal{S}^{\mathcal{M}}_{\text{Roum}}}$ respectively for the Beurling-case. If $(s_k)_k\in\hyperlink{S-beur}{\tilde{\mathcal{S}}^{\mathcal{M}}_{\text{Beur}}},\hyperlink{s-beur}{\mathcal{S}^{\mathcal{M}}_{\text{Beur}}}$ holds for $x\in\Lambda$, then also for all $y\le x$, too. All occurring sets are stable w.r.t. $(\cdot_k)_k\mapsto(B^k\cdot_k)_k$ for arbitrary $B>0$.

Using \cite[Lemma 4.6.]{KMRplot} directly we get:

\begin{proposition}\label{existencequantor5}
Let $\mathcal{M}=\{M^x: x\in\Lambda\}$ be \hyperlink{Marb}{$(\mathcal{M})$}, $E,F$ be Banach spaces, $U\subseteq E$ open and $f:U\rightarrow F$ a $\mathcal{E}$-mapping. Then the following are equivalent:
\begin{itemize}
\item[$(1)$] $f$ is $\mathcal{E}_{\{\mathcal{M}\}}=\mathcal{E}_{\{\mathcal{M}\}}^{\on{b}}$.

\item[$(2)$] For each compact $K\subseteq U$ there exists $x\in\Lambda$ such that for each $(r_k)_k\in\hyperlink{r-roum}{\mathcal{R}_{\on{Roum}}}$
$$\left\{\frac{f^{(k)}(a)(v_1,\dots,v_k) r_k}{M^x_k}: a\in K, k\in\NN, \|v_i\|_E\le 1\right\}$$
is bounded in $F$.

\item[$(3)$] For each compact $K\subseteq U$ there exists $x\in\Lambda$ such that for each $(r_k)_k\in\hyperlink{r-roum-sub}{\mathcal{R}_{\on{Roum},\on{sub}}}$ there exists $\varepsilon>0$ such that
$$\left\{\frac{f^{(k)}(a)(v_1,\dots,v_k) r_k\varepsilon^k}{M^x_k}: a\in K, k\in\NN, \|v_i\|_E\le 1\right\}$$
is bounded in $F$.
\end{itemize}
\end{proposition}
Note that $\mathcal{E}_{\{\mathcal{M}\}}=\mathcal{E}_{\{\mathcal{M}\}}^{\on{b}}$ holds by Lemma \ref{Plots-existencequantifiers} below, but for our approach in this work we also have to get rid of the second existence quantifier.

\subsection{Roumieu-case with $(\mathcal{M}_{\{\on{FdB}\}})$}\label{subsection31}
We prove the following generalization of \cite[Lemma 9.2.]{KM97}:
\begin{lemma}\label{existencequantor1}
Let $\mathcal{M}=\{M^x: x\in\Lambda\}$ be \hyperlink{Msc}{$(\mathcal{M}_{\on{sc}})$} with $\Lambda=\NN_{>0}$ and \hyperlink{R-FdB}{$(\mathcal{M}_{\{\on{FdB}\}})$}. For a formal power series $\sum_{k\ge 0}a^x_k t^k=\sum_{k\ge 0}\frac{b_k}{k! m^x_k} t^k$, so $a^x_k:=\frac{b_k}{M^x_k}$, the following are equivalent:

\begin{itemize}
\item[$(1)$] There exists $x\in\Lambda$ such that $\sum_{k\ge 0}a^x_k t^k$ has positive radius of convergence.

\item[$(2)$] $\sum_{k\ge 0}\frac{b_k r_k s_k}{k!}$ converges absolutely for all $(r_k)_k\in\hyperlink{r-roum}{\mathcal{R}_{\on{Roum}}}$ and $(s_k)_k\in\hyperlink{s-roum}{\mathcal{S}^{\mathcal{M}}_{\on{Roum}}}$.

\item[$(3)$] The sequence $\left(\frac{b_k r_k s_k}{k!}\right)_k$ is bounded for all $(r_k)_k\in\hyperlink{r-roum}{\mathcal{R}_{\on{Roum}}}$ and $(s_k)_k\in\hyperlink{s-roum}{\mathcal{S}^{\mathcal{M}}_{\on{Roum}}}$.

\item[$(4)$] For each $(r_k)_k\in\hyperlink{r-roum-sub}{\mathcal{R}_{\on{Roum},\on{sub}}}$ and for each $(s_k)_k\in\hyperlink{s-roum-FdB}{\mathcal{S}^{\mathcal{M}}_{\on{Roum},\on{FdB}}}$ there exists $\varepsilon>0$ such that $\left(\frac{b_k r_k s_k}{k!}\varepsilon^k\right)_k$ is bounded.
\end{itemize}
\end{lemma}

\demo{Proof}
$(1)\Rightarrow(2)$ For the given series ($x\in\Lambda$ coming from $(1)$) and arbitrary $(r_k)_k$ and $(s_k)_k$ as considered in $(2)$ we have
$$\sum_{k\ge 0}\frac{b_k r_k s_k}{k!}=\sum_{k\ge 0}a^x_k m^x_k r_k s_k=\sum_{k\ge 0}(a^x_k t^k)\underbrace{(s_k m^x_k)}_{\le C_x^k}\frac{r_k}{t^k}\le\sum_{k\ge 0}(a^x_k t^k)\underbrace{r_k\left(\frac{C_x}{t}\right)^k}_{\rightarrow 0,\text{as}\;k\rightarrow\infty},$$
hence the first sum converges for $t>0$ sufficiently small.

$(2)\Rightarrow(3)\Rightarrow(4)$ are clearly satisfied.

$(4)\Rightarrow(1)$ Since \hyperlink{R-FdB}{$(\mathcal{M}_{\{\text{FdB}\}})$} is satisfied and $m^{x}\le m^{y}$ for $x\le y$ we can associate to each $x\in\Lambda$ the index $\alpha(x):=\min\{y\in\Lambda: (m^x)^{\circ}\hyperlink{mpreceq}{\preceq}m^y\}$. Since $(m^{x})^{\circ}\le (m^{y})^{\circ}$ for $x\le y$ we also have $\alpha(x)\le\alpha(y)$ for such indices and $\lim_{x\rightarrow\infty}\alpha(x)=+\infty$.

On the other hand for $y\ge\alpha(1)$ we can define $\beta(y):=\max\{x\in\Lambda:\alpha(x)\le y\}$ which is clearly well-defined. So $\beta(y_1)\le\beta(y_2)$ for $y_1\le y_2$, $\lim_{y\rightarrow\infty}\beta(y)=+\infty$ and finally by construction for each $x\in\NN_{>0}$, $x\ge\alpha(1)$, there exist $y\in\NN_{>0}$, $y\le x$, with $(m^y)^{\circ}\hyperlink{mpreceq}{\preceq}m^x$. Note that this does not imply \hyperlink{B-FdB}{$(\mathcal{M}_{(\text{FdB})})$}. W.l.o.g. we could assume that $\alpha(x)=x+1$ and so $\beta(y)=y-1$. If $\mathcal{M}$ has in addition \hyperlink{R-Comega}{$(\mathcal{M}_{\{\mathcal{C}^{\omega}\}})$}, i.e. the real analytic functions are contained in $\mathcal{E}_{\{\mathcal{M}\}}$, then we can take w.l.o.g. $M^1=(p!)_{p\in\NN}$, so $m^1_p=1$ for each $p$ and $\alpha(1)=1$.\vspace{6pt}

We prove by contradiction. So assume that each $\sum_{k\ge 0}a^x_k t^k$ would have radius of convergence $0$. Then we would get $\sum_{k\ge 0}|a^x_k|\left(\frac{1}{n^2}\right)^k=+\infty$ for each $n\in\NN_{>0}$ and each $x\in\Lambda=\NN_{>0}$. Consider now $n\in\NN_{>0}$ and $x:=n+\alpha(1)$ and so we find an increasing sequence $(k_n)_{n\ge 0}$ with $k_0=1$, $\lim_{n\rightarrow\infty}k_n=+\infty$ such that
\begin{equation}\label{FdBproj}
\forall\;n\in\NN_{>0}:\;\;\sum_{k=k_{n-1}}^{k_{n}-1}|a^{n+\alpha(1)}_k|\left(\frac{1}{n^2}\right)^k\ge 1.
\end{equation}
We put now
$$r_k:=\left(\frac{1}{n^2}\right)^k\;\;\text{for}\;\;k_{n-1}\le k\le k_n-1, n\in\NN_{>0},$$
and show $(r_k)_k\in\hyperlink{r-roum-sub}{\mathcal{R}_{\text{Roum},\text{sub}}}$. For $k_{n-1}\le k\le k_n-1$ by definition $r_k t^k=\left(\frac{t}{n^2}\right)^k$, and so $r_k t^k\rightarrow 0$ as $k\rightarrow\infty$ and all $t>0$. Clearly $(r_k)_k$ is also log-sub-additive. In addition one can see that $(\sqrt{r_k})_k\in\hyperlink{r-roum-sub}{\mathcal{R}_{\text{Roum},\text{sub}}}$ and so for all $\varepsilon>0$ there exists $k_{\varepsilon}\in\NN$ such that for all $k\ge k_{\varepsilon}$ we have $\sqrt{r_k}\frac{1}{\varepsilon^k}\le 1\Leftrightarrow\sqrt{r_k}\le\varepsilon^k$.\vspace{6pt}

No we define $s:=(s_k)_k$. We put $s_k:=\frac{1}{m^{\gamma(k)}_{k}}$, where $\gamma(k):=n+\alpha(1)$ for $k_{n-1}\le k\le k_n-1$, $n\in\NN_{>0}$, and show $(s_k)_k\in\hyperlink{s-roum-FdB}{\mathcal{S}^{\mathcal{M}}_{\text{Roum},\text{FdB}}}$.

So let $x\in\Lambda$ be arbitrary (large) but fixed, then for $k_{n-1}\le k\le k_n-1$ we get $s_k m^x_k=\frac{m^x_k}{m^{\gamma(k)}_k}=\frac{M^x_k}{M^{\gamma(k)}_k}$. For all $k\in\NN$ we can estimate $\frac{M^x_k}{M^{\gamma(k)}_k}\le C_x^k$ with some constant $C_x>0$, because $\lim_{k\rightarrow\infty}\gamma(k)=+\infty$. This proves $(s_k)_k\in\hyperlink{s-roum}{\mathcal{S}^{\mathcal{M}}_{\text{Roum}}}$. Define
$$\hat{s}_k:=\frac{1}{m^{\beta(\gamma(k))}_k}\hspace{20pt}\text{for}\;\;k_{n-1}\le k\le k_n-1, n\in\NN_{>0},$$
and similarly we find a constant $D_x>0$ such that $\hat{s}_k m^x_k\le D^k_x$ for each $x\in\Lambda$ and $k\in\NN$ because $\lim_{k\rightarrow\infty}\beta(\gamma(k))=+\infty$. This proves $(\hat{s}_k)_k\in\hyperlink{s-roum}{\mathcal{S}^{\mathcal{M}}_{\text{Roum}}}$. For $\delta_1+\dots+\delta_j=k$ we obtain for $k\in\NN$ with $k_{n-1}\le k\le k_n-1$, $n\in\NN_{>0}$:
\begin{align*}
s_k&=\frac{1}{m^{\gamma(k)}_k}\le C h^k\frac{1}{m_j^{\beta(\gamma(k))} m_{\delta_1}^{\beta(\gamma(k))}\cdots m_{\delta_j}^{\beta(\gamma(k))}}
\\&
\le C h^k\frac{1}{m_j^{\beta(\gamma(j))} m_{\delta_1}^{\beta(\gamma(\delta_1))}\cdots m_{\delta_j}^{\beta(\gamma(\delta_j))}}=C h^k\hat{s}_j\hat{s}_{\delta_1}\cdots\hat{s}_{\delta_j},
\end{align*}
which precisely shows $s\hyperlink{mpreceq}{\preceq}\hat{s}_{o}$. The first inequality holds by \hyperlink{R-FdB}{$(\mathcal{M}_{\{\text{FdB}\}})$} and by definition of $\beta$, the second because $j,\delta_1,\dots,\delta_j\le k$. So $s$ is as desired.

Moreover
$$\sum_{k\ge 1}\frac{|b_k| r_k s_k}{k!}=\sum_{n\ge 1}\sum_{k=k_{n-1}}^{k_n-1}\frac{|b_k| r_k s_k}{k!}=\sum_{n\ge 1}\sum_{k=k_{n-1}}^{k_{n}-1}|a^{n+\alpha(1)}_k|\left(\frac{1}{n^2}\right)^k\ge\sum_{n\ge 1}1=+\infty,$$
because by definition $\frac{|b_k|}{k!} s_k=\frac{|b_k|}{M^{n+\alpha(1)}_k}=|a^{n+\alpha(1)}_k|$ for $k_{n-1}\le k\le k_n-1$ (note that $n(k)=n+\alpha(1)$ for $k\in[k_{n-1},k_n-1]$).\vspace{6pt}

Finally we show that $\left(\frac{b_k}{k!}\sqrt{r_k} s_k(2\varepsilon)^k\right)_k$ cannot be bounded for any $\varepsilon>0$. First we get
\begin{align*}
\sum_{k\ge 1}\frac{|b_k|}{k!}\sqrt{r_k} s_k\varepsilon^k\ge\sum_{k\ge k_{\varepsilon}}\frac{|b_k|}{k!}\sqrt{r_k} s_k\underbrace{\varepsilon^k}_{\ge\sqrt{r_k}}\ge\sum_{k\ge k_{\varepsilon}}\frac{|b_k|}{k!} r_k s_k=+\infty.
\end{align*}
But if the sequence would be bounded for some $\varepsilon$, then for all $k\in\NN$ we would get $\frac{b_k}{k!}\sqrt{r_k} s_k\varepsilon^k\le\frac{C}{2^k}$, hence $\sum_{k\ge 0}\frac{|b_k|}{k!}\sqrt{r_k} s_k\varepsilon^k\le\sum_{k\ge 0}\frac{C}{2^k}=2C$, a contradiction.
\qed\enddemo

We use Lemma \ref{existencequantor1} to generalize \cite[Lemma 4.6.]{KMRplot}.

\begin{proposition}\label{existencequantorroumieustrong}
Let $\mathcal{M}=\{M^x: x\in\Lambda\}$ be \hyperlink{Msc}{$(\mathcal{M}_{\on{sc}})$} with $\Lambda=\NN_{>0}$ and \hyperlink{R-FdB}{$(\mathcal{M}_{\{\on{FdB}\}})$}. Let $E,F$ be Banach spaces, $U\subseteq E$ open and $f:U\rightarrow F$ a $\mathcal{E}$-mapping. Then the following are equivalent:
\begin{itemize}
\item[$(1)$] $f$ is $\mathcal{E}_{\{\mathcal{M}\}}=\mathcal{E}_{\{\mathcal{M}\}}^{\on{b}}$.

\item[$(2)$] For each compact $K\subseteq U$, for each $(r_k)_k\in\hyperlink{r-roum}{\mathcal{R}_{\on{Roum}}}$ and each $(s_k)_k\in\hyperlink{s-roum}{\mathcal{S}^{\mathcal{M}}_{\on{Roum}}}$ the set
$$\left\{\frac{f^{(k)}(a)(v_1,\dots,v_k)}{k!} r_k s_k: a\in K, k\in\NN, \|v_i\|_E\le 1\right\}$$
is bounded in $F$.

\item[$(3)$] For each compact $K\subseteq U$, for each $(r_k)_k\in\hyperlink{r-roum-sub}{\mathcal{R}_{\on{Roum},\on{sub}}}$ and for each $(s_k)_k\in\hyperlink{s-roum-FdB}{\mathcal{S}^{\mathcal{M}}_{\on{Roum},\on{FdB}}}$, there exists $\varepsilon>0$ such that the set
$$\left\{\frac{f^{(k)}(a)(v_1,\dots,v_k)}{k!} r_k s_k\varepsilon^k: a\in K, k\in\NN, \|v_i\|_E\le 1\right\}$$
is bounded in $F$.
\end{itemize}
\end{proposition}

\demo{Proof}
$(1)\Rightarrow(2)$ Let $f$ be $\mathcal{E}_{\{\mathcal{M}\}}$ and $K\subseteq U$ compact, then estimate as follows (where we use Lemma \ref{Plots-existencequantifiers} below):
\begin{align*}
\left\|\frac{f^{(k)}(a)}{k!}r_k s_k\right\|_{L^k(E,F)}&=\left\|\frac{f^{(k)}(a)}{k! m^x_k h^k}\right\|_{L^k(E,F)}|r_k h^k\underbrace{s_k m^x_k}_{\le C_x^k}|\le\left\|\frac{f^{(k)}(a)}{h^k M^x_k}\right\|_{L^k(E,F)}\underbrace{\left|r_k(C_x h)^k\right|}_{\rightarrow 0,\;\text{as}\;k\rightarrow\infty}
\end{align*}
for $a\in K$, $x\in\Lambda$ and $h>0$ large enough (depending on $K$ and $f$) and for arbitrary $(r_k)_k$ and $(s_k)_k$ as considered in $(2)$.

$(2)\Rightarrow(3)$ Take $\varepsilon=1$.

$(3)\Rightarrow(1)$ We use $(4)\Rightarrow(1)$ in Lemma \ref{existencequantor1}. Let $K\subseteq U$ be an arbitrary compact set but fixed and put $b_k:=\sup_{a\in K}\left\|f^{(k)}(a)\right\|_{L^k(E,F)}$. Then there exists $h>0$ and $x\in\Lambda$ such that $\sup_{k\in\NN}\frac{b_k}{M_x^kh^k}<+\infty$, hence $f$ is $\mathcal{E}_{\{\mathcal{M}\}}$.
\qed\enddemo

\subsection{Roumieu-case without $(\mathcal{M}_{\{\on{FdB}\}})$}\label{subsection32}
\begin{lemma}\label{existencequantor2}
Let $\mathcal{M}=\{M^x: x\in\Lambda\}$ be \hyperlink{Msc}{$(\mathcal{M}_{\on{sc}})$} with $\Lambda=\NN_{>0}$. For a formal power series $\sum_{k\ge 0}a^x_k t^k=\sum_{k\ge 0}\frac{b_k}{M^x_k}t^k$ the following are equivalent:

\begin{itemize}
\item[$(1)$] There exists $x\in\Lambda$ such that $\sum_{k\ge 0}a^x_k t^k$ has positive radius of convergence.

\item[$(2)$] $\sum_{k\ge 0}b_k r_k s_k$ converges absolutely for all $(r_k)_k\in\hyperlink{r-roum}{\mathcal{R}_{\on{Roum}}}$ and $(s_k)_k\in\hyperlink{S-roum}{\tilde{\mathcal{S}}^{\mathcal{M}}_{\on{Roum}}}$.

\item[$(3)$] $(b_k r_k s_k)_k$ is bounded for all $(r_k)_k\in\hyperlink{r-roum}{\mathcal{R}_{\on{Roum}}}$ and $(s_k)_k\in\hyperlink{S-roum-sub}{\tilde{\mathcal{S}}^{\mathcal{M}}_{\on{Roum},\on{sub}}}$.

\item[$(4)$] For each $(r_k)_k\in\hyperlink{r-roum-sub}{\mathcal{R}_{\on{Roum},\on{sub}}}$ and for each $(s_k)_k\in\hyperlink{S-roum-sub}{\tilde{\mathcal{S}}^{\mathcal{M}}_{\on{Roum},\on{sub}}}$ there exists $\varepsilon>0$ such that $(b_k r_k s_k\varepsilon^k)_k$ is bounded.
\end{itemize}
If $\mathcal{M}$ is \hyperlink{Marb}{$(\mathcal{M})$}, then in $(3)$ and $(4)$ we replace \hyperlink{S-roum-sub}{$\tilde{\mathcal{S}}^{\mathcal{M}}_{\on{Roum},\on{sub}}$} by \hyperlink{S-roum}{$\tilde{\mathcal{S}}^{\mathcal{M}}_{\on{Roum}}$}.
\end{lemma}

\demo{Proof}
$(1)\Rightarrow(2)\Rightarrow(3)\Rightarrow(4)$ is the same as in Lemma \ref{existencequantor1}. For $(4)\Rightarrow(1)$ we prove again by contradiction. In \eqref{FdBproj} consider $x=n\in\NN_{>0}$, take the same $r=(r_k)_k$ and for $s=(s_k)_k$ we put $s_k:=\frac{1}{M^n_k}$ if $k_{n-1}\le k\le k_n-1$. If $\mathcal{M}$ is \hyperlink{Msc}{$(\mathcal{M}_{\on{sc}})$} then we have $M^x_j M^x_k\le M^x_{j+k}$ for each $j,k\in\NN$ and $x\in\Lambda$ and $M^x\le M^y$ for $x\le y$. This implies $(s_k)_k\in\hyperlink{S-roum-sub}{\tilde{\mathcal{S}}^{\mathcal{M}}_{\on{Roum},\on{sub}}}$. If $\mathcal{M}$ is \hyperlink{Marb}{$(\mathcal{M})$}, then $(s_k)_k\in\hyperlink{S-roum}{\tilde{\mathcal{S}}^{\mathcal{M}}_{\on{Roum}}}$ holds by definition.
\qed\enddemo

So we can prove a new version of Proposition \ref{existencequantorroumieustrong}.

\begin{proposition}\label{existencequantorroumieuweak}
Let $\mathcal{M}=\{M^x: x\in\Lambda\}$ be \hyperlink{Msc}{$(\mathcal{M}_{\on{sc}})$} with $\Lambda=\NN_{>0}$. Let $E,F$ be Banach spaces, $U\subseteq E$ open and $f:U\rightarrow F$ a $\mathcal{E}$-mapping, then the following are equivalent:

\begin{itemize}
\item[$(1)$] $f$ is $\mathcal{E}_{\{\mathcal{M}\}}=\mathcal{E}_{\{\mathcal{M}\}}^{\on{b}}$.

\item[$(2)$] For each compact $K\subseteq U$, for each $(r_k)_k\in\hyperlink{r-roum}{\mathcal{R}_{\on{Roum}}}$ and for each $(s_k)_k\in\hyperlink{S-roum}{\tilde{\mathcal{S}}^{\mathcal{M}}_{\on{Roum}}}$ the set
 $$\left\{f^{(k)}(a)(v_1,\dots,v_k) r_k s_k: a\in K, k\in\NN, \|v_i\|_E\le 1\right\}$$
is bounded in $F$.

\item[$(3)$] For each compact $K\subseteq U$, for each $(r_k)_k\in\hyperlink{r-roum-sub}{\mathcal{R}_{\on{Roum},\on{sub}}}$ and for each $(s_k)_k\in\hyperlink{S-roum-sub}{\tilde{\mathcal{S}}^{\mathcal{M}}_{\on{Roum},\on{sub}}}$, there exists $\varepsilon>0$ such that the set
    $$\left\{f^{(k)}(a)(v_1,\dots,v_k) r_k s_k\varepsilon^k: a\in K, k\in\NN, \|v_i\|_E\le 1\right\}$$
is bounded in $F$.
\end{itemize}
If $\mathcal{M}$ is \hyperlink{Marb}{$(\mathcal{M})$}, then in $(3)$ we replace \hyperlink{S-roum-sub}{$\tilde{\mathcal{S}}^{\mathcal{M}}_{\on{Roum},\on{sub}}$} by \hyperlink{S-roum}{$\tilde{\mathcal{S}}^{\mathcal{M}}_{\on{Roum}}$}.
\end{proposition}

\demo{Proof}
Use precisely the same arguments as in Proposition \ref{existencequantorroumieustrong}, for $(3)\Rightarrow(1)$ we use $(4)\Rightarrow(1)$ in Lemma \ref{existencequantor2}.
\qed\enddemo

\subsection{Beurling-case with $(\mathcal{M}_{(\on{FdB})})$}\label{subsection33}
\begin{lemma}\label{existencequantor3}
Let $\mathcal{M}=\{M^x: x\in\Lambda\}$ be \hyperlink{Msc}{$(\mathcal{M}_{\on{sc}})$} with $\Lambda=\RR_{>0}$ and \hyperlink{B-FdB}{$(\mathcal{M}_{(\on{FdB})})$}. For a formal power series $\sum_{k\ge 0}a^x_k t^k=\sum_{k\ge 0}\frac{b_k}{M^x_k} t^k$, so $a^x_k:=\frac{b_k}{M^x_k}$, the following are equivalent:

\begin{itemize}
\item[$(1)$] The series $\sum_{k\ge 0}a^x_k t^k$ has infinite radius of convergence for each $x\in\Lambda$.

\item[$(2)$] For each $(r_k)_k\in\hyperlink{r-beur-sub}{\mathcal{R}_{\on{Beur},\on{sub}}}$ and for each $(s_k)_k\in\hyperlink{s-beur-FdB}{\mathcal{S}^{\mathcal{M}}_{\on{Beur},\on{FdB}}}$ the sequence $\left(\frac{b_k}{k!} r_k s_k\delta^k\right)_k$ is bounded for each $\delta>0$.
\end{itemize}
\end{lemma}

\demo{Proof}
$(1)\Rightarrow(2)$ Let $(r_k)_k$ and $(s_k)_k$ be given as considered in $(2)$, then
$$\sum_{k\ge 0}\frac{b_k}{k!} r_k s_k\delta^k=\sum_{k\ge 0} a_k^x\underbrace{(m^x_k s_k)}_{\le C_x^k}(r_k t^k)\left(\frac{\delta}{t}\right)^k\le\sum_{k\ge 0} a_k^x\underbrace{(r_k t^k)}_{\rightarrow 0,\;\text{as}\;k\rightarrow\infty}\left(\frac{\delta C_x}{t}\right)^k$$
is absolutely convergent for each $\delta>0$. The index $x\in\Lambda$ was chosen such that $s_k m^x_k\le C_x^k$ holds for all $k\in\NN$ and it is depending on $(s_k)_k\in\hyperlink{s-beur}{\mathcal{S}^{\mathcal{M}}_{\on{Beur}}}$. The real number $t>0$ was chosen in such a way that $r_k t^k\rightarrow 0$ as $k\rightarrow\infty$. Hence $\left(\frac{b_k}{k!} r_k s_k\delta^k\right)_k$ is bounded for each $\delta>0$.\vspace{6pt}

$(2)\Rightarrow(1)$ Assume that there would exist $x\in\Lambda$ such that $\sum_{k\ge 0}a^x_k t^k$ would have finite radius of convergence. Then there would exist $h>0$ such that $\sum_{k\ge 0}|a^x_k| n^k=+\infty$ for each $n>h$. Put now $r_k:=\frac{1}{n^k}$ for some $n>h$ and $s_k:=\frac{1}{m^x_k}$.

Clearly $(r_k)_k\in\hyperlink{r-beur-sub}{\mathcal{R}_{\text{Beur},\text{sub}}}$ holds.

Also $(s_k)_k$ is as desired. By \hyperlink{B-FdB}{$(\mathcal{M}_{(\text{FdB})})$} for all $x\in\Lambda$ there exists $y\in\Lambda$ and $D>0$ such that for all $\alpha_1+\dots+\alpha_j=k$ we get
$$s_k:=\frac{1}{m^x_k}\le D^k\frac{1}{m^{y}_j m^{y}_{\alpha_1}\cdots m^{y}_{\alpha_j}}=:D^k\hat{s}_j\hat{s}_{\alpha_1}\cdots\hat{s}_{\alpha_j},$$
where we have put $\hat{s}_j:=\frac{1}{m^{y}_j}$. We have $y\le x$, since $(m^{y})^{\circ}\le(m^{x})^{\circ}$ for $y\le x$. Clearly $(\hat{s}_k)_k\in\hyperlink{s-beur}{\mathcal{S}^{\mathcal{M}}_{\text{Beur}}}$, hence $(s_k)_k\in\hyperlink{s-beur-FdB}{\mathcal{S}^{\mathcal{M}}_{\text{Beur},\text{FdB}}}$ and so both sequences are as considered in $(2)$. But then there would exist $C>0$ such that for all $k\in\NN$:
$$C>\frac{b_k}{k!} s_k r_k(2n^2)^k=\frac{b_k}{k! m^x_k} r_k(2n^2)^k=a^x_k r_k n^{2k} 2^k=a_k^x n^k 2^k.$$
Hence $\sum_{k\ge 0}|a^x_k| n^k\le C\sum_{k\ge 0}\frac{1}{2^k}=2C$, a contradiction.
\qed\enddemo

Using the previous result we can show:

\begin{proposition}\label{existencequantorbeurlingstrong}
Let $\mathcal{M}=\{M^x: x\in\Lambda\}$ be \hyperlink{Msc}{$(\mathcal{M}_{\on{sc}})$} with $\Lambda=\RR_{>0}$ and \hyperlink{B-FdB}{$(\mathcal{M}_{(\on{FdB})})$}. Let $E,F$ be Banach spaces, $U\subseteq E$ open and $f:U\rightarrow F$ a $\mathcal{E}$-mapping, then the following are equivalent:
\begin{itemize}
\item[$(1)$] $f$ is $\mathcal{E}_{(\mathcal{M})}=\mathcal{E}_{(\mathcal{M})}^{\on{b}}$.

\item[$(2)$] For each compact $K\subseteq U$, for each $(r_k)_k\in\hyperlink{r-beur}{\mathcal{R}_{\on{Beur}}}$ and for each $(s_k)_k\in\hyperlink{s-beur}{\mathcal{S}^{\mathcal{M}}_{\on{Beur}}}$ the set
$$\left\{\frac{f^{(k)}(a)(v_1,\dots,v_k)}{k!} r_k s_k: a\in K, k\in\NN, \|v_i\|_E\le 1\right\}$$
is bounded in $F$.

\item[$(3)$] For each compact $K\subseteq U$, for each $(r_k)_k\in\hyperlink{r-beur-sub}{\mathcal{R}_{\on{Beur},\on{sub}}}$ and for each $(s_k)_k\in\hyperlink{s-beur-FdB}{\mathcal{S}^{\mathcal{M}}_{\on{Beur},\on{FdB}}}$ the set
$$\left\{\frac{f^{(k)}(a)(v_1,\dots,v_k)}{k!} r_k s_k\delta^k: a\in K, k\in\NN, \|v_i\|_E\le 1\right\}$$
is bounded in $F$ for each $\delta>0$.
\end{itemize}
\end{proposition}

\demo{Proof}
$(1)\Rightarrow(2)$ Let $f$ be $\mathcal{E}_{(\mathcal{M})}$ and $(r_k)_k$, $(s_k)_k$ given by $(2)$, then we can estimate as follows (where we use Lemma \ref{Plots-allquantifiers} below):
\begin{align*}
\left\|\frac{f^{(k)}(a)}{k!} r_k s_k\right\|_{L^k(E,F)}&=\left\|\frac{f^{(k)}(a)}{k! m^x_k h^k}\right\|_{L^k(E,F)}|r_k h^k\underbrace{s_k m^x_k}_{\le C_x^k}|\le\left\|\frac{f^{(k)}(a)}{M^x_k h^k}\right\|_{L^k(E,F)}\underbrace{\left|r_k(C_x h)^k\right|}_{\rightarrow 0,\;\text{as}\;k\rightarrow\infty}
\end{align*}
for $a\in K$. We have chosen $x\in\Lambda$ depending on $(s_k)_k\in\hyperlink{s-beur}{\mathcal{S}^{\mathcal{M}}_{\on{Beur}}}$ such that $s_k m^x_k\le C_x^k$ and $h>0$ depending on given $(r_k)_k\in\hyperlink{r-beur}{\mathcal{R}_{\on{Beur}}}$ such that $r_k(C_xh)^k\rightarrow 0$ as $k\rightarrow\infty$.

$(2)\Rightarrow(3)$ Replace in $(2)$ the sequence $(r_k)_k$ by $(r_k\delta^k)_k$.

$(3)\Rightarrow(1)$ Use $(2)\Rightarrow(1)$ in Lemma \ref{existencequantor3}. Let $K\subseteq U$ be a compact set, arbitrary but fixed. Then put $b_k:=\sup_{a\in K}\left\|f^{(k)}(a)\right\|_{L^k(E,F)}$ and so for each $h>0$ and each $x\in\Lambda$ we have that $\sup_{k\in\NN}\frac{b_k}{M_x^k h^k}<+\infty$, hence $f$ is $\mathcal{E}_{(\mathcal{M})}$.
\qed\enddemo

\subsection{Beurling-case without $(\mathcal{M}_{(\on{FdB})})$}\label{subsection34}
\begin{lemma}\label{existencequantor4}
Let $\mathcal{M}=\{M^x: x\in\Lambda\}$ be \hyperlink{Msc}{$(\mathcal{M}_{\on{sc}})$} with $\Lambda=\RR_{>0}$. For a formal power series $\sum_{k\ge 0}a^x_k t^k=\sum_{k\ge 0}\frac{b_k}{M^x_k} t^k$, $a^x_k:=\frac{b_k}{M^x_k}$, the following are equivalent:

\begin{itemize}
\item[$(1)$] The series $\sum_{k\ge 0}a^x_k t^k$ has infinite radius of convergence for each $x\in\Lambda$.

\item[$(2)$] For each $(r_k)_k\in\hyperlink{r-beur-sub}{\mathcal{R}_{\on{Beur},\on{sub}}}$ and for each $(s_k)_k\in\hyperlink{S-beur-sub}{\tilde{\mathcal{S}}^{\mathcal{M}}_{\on{Beur},\on{sub}}}$ the sequence $(b_k r_k s_k\delta^k)_k$ is bounded for each $\delta>0$.
\end{itemize}
If $\mathcal{M}$ is \hyperlink{Marb}{$(\mathcal{M})$}, then in $(2)$ we replace \hyperlink{S-beur-sub}{$\tilde{\mathcal{S}}^{\mathcal{M}}_{\on{Beur},\on{sub}}$} by \hyperlink{S-beur}{$\tilde{\mathcal{S}}^{\mathcal{M}}_{\on{Beur}}$}.
\end{lemma}

\demo{Proof}
Proceed as in Lemma \ref{existencequantor3}: For $(2)\Rightarrow(1)$ we put $s_k:=\frac{1}{M^x_k}$, where $x\in\Lambda$ is the index arising by the contradiction argument.

Hence $(s_k)_k\in\hyperlink{S-beur-sub}{\tilde{\mathcal{S}}^{\mathcal{M}}_{\text{beur},\text{sub}}}$ holds whenever $\mathcal{M}$ is \hyperlink{Msc}{$(\mathcal{M}_{\on{sc}})$} since each $M^x$ is log-convex. If $\mathcal{M}$ is \hyperlink{Marb}{$(\mathcal{M})$}, then $(s_k)_k\in\hyperlink{S-beur}{\tilde{\mathcal{S}}^{\mathcal{M}}_{\on{Beur}}}$ is clear.
\qed\enddemo

So we are able to prove:

\begin{proposition}\label{existencequantorbeurlingweak}
Let $\mathcal{M}=\{M^x: x\in\Lambda\}$ be \hyperlink{Msc}{$(\mathcal{M}_{\on{sc}})$} with $\Lambda=\RR_{>0}$. Let $E,F$ be Banach spaces, $U\subseteq E$ open and $f:U\rightarrow F$ a $\mathcal{E}$-mapping, then the following are equivalent:

\begin{itemize}
\item[$(1)$] $f$ is $\mathcal{E}_{(\mathcal{M})}=\mathcal{E}_{(\mathcal{M})}^{\on{b}}$.

\item[$(2)$] For each compact $K\subseteq U$, for each $(r_k)_k\in\hyperlink{r-beur}{\mathcal{R}_{\on{Beur}}}$ and for each $(s_k)_k\in\hyperlink{S-beur}{\tilde{\mathcal{S}}^{\mathcal{M}}_{\on{Beur}}}$ the set
$$\left\{f^{(k)}(a)(v_1,\dots,v_k) r_k s_k: a\in K, k\in\NN, \|v_i\|_E\le 1\right\}$$
is bounded in $F$.

\item[$(3)$] For each compact $K\subseteq U$, for each $(r_k)_k\in\hyperlink{r-beur-sub}{\mathcal{R}_{\on{Beur},\on{sub}}}$ and for each $(s_k)_k\in\hyperlink{S-beur-sub}{\tilde{\mathcal{S}}^{\mathcal{M}}_{\on{Beur},\on{sub}}}$ the set
$$\left\{f^{(k)}(a)(v_1,\dots,v_k) r_k s_k\delta^k: a\in K, k\in\NN, \|v_i\|_E\le 1\right\}$$
is bounded in $F$ for each $\delta>0$.
\end{itemize}
If $\mathcal{M}$ is \hyperlink{Marb}{$(\mathcal{M})$}, then in $(3)$ we replace \hyperlink{S-beur-sub}{$\tilde{\mathcal{S}}^{\mathcal{M}}_{\on{Beur},\on{sub}}$} by \hyperlink{S-beur}{$\tilde{\mathcal{S}}^{\mathcal{M}}_{\on{Beur}}$}.
\end{proposition}

\demo{Proof}
The proof is the same as for Proposition \ref{existencequantorbeurlingstrong}. For $(3)\Rightarrow(1)$ we use $(2)\Rightarrow(1)$ in Lemma \ref{existencequantor4}.
\qed\enddemo

\section{Closedness under composition}\label{section4}

\subsection{First observations}\label{subsection41}

First we generalize \cite[Lemma 4.2.]{KMRplot}:

\begin{lemma}\label{Plots-allquantifiers}
Let $\mathcal{M}$ be \hyperlink{Marb}{$(\mathcal{M})$}, then $\mathcal{E}_{(\mathcal{M})}=\mathcal{E}_{(\mathcal{M})}^{\on{b}}$.
\end{lemma}

\demo{Proof}
Let $E,F$ be convenient, $U\subseteq E$ a $c^{\infty}$-open subset and let $f:U\rightarrow F$ be a $\mathcal{E}$-mapping. Then we obtain the following equivalences, where the set $B$ runs through all closed absolutely convex bounded subsets in $E$ and $K$ runs through all sets in $U_B$ which are compact w.r.t. the norm $\|\cdot\|_B$:
\begin{align*}
&f\in\mathcal{E}_{(\mathcal{M})}(U,F)
\\&
\Longleftrightarrow\forall\;\alpha\in F^{\ast}\;\forall\;B\;\forall\;K\subseteq U_B\;\forall\;x\in\Lambda\;\forall\;h>0:
\\&
\left\{\frac{(\alpha\circ f)^{(k)}(a)(v_1,\dots,v_k)}{h^k M^x_k}: a\in K, k\in\NN, \|v_i\|_B\le 1\right\}\;\text{is bounded in}\;\RR
\\&
\Longleftrightarrow\forall\;B\;\forall\;K\subseteq U_B\;\forall\;x\in\Lambda\;\forall\;h>0\;\forall\;\alpha\in F^{\ast}:
\\&
\alpha\left(\left\{\frac{f^{(k)}(a)(v_1,\dots,v_k)}{h^k M^x_k}: a\in K, k\in\NN, \|v_i\|_B\le 1\right\}\right)\;\text{is bounded in}\;\RR
\\&
\Longleftrightarrow\forall\;B\;\forall\;K\subseteq U_B\;\forall\;x\in\Lambda\;\forall\;h>0:
\\&
\left\{\frac{f^{(k)}(a)(v_1,\dots,v_k)}{h^k M^x_k}: a\in K, k\in\NN, \|v_i\|_B\le 1\right\}\;\text{is bounded in}\;\RR
\\&
\Longleftrightarrow f\in\mathcal{E}_{(\mathcal{M})}^{\on{b}}(U,F).
\end{align*}
\qed\enddemo

But in general we do not have $\mathcal{E}_{\{\mathcal{M}\}}=\mathcal{E}_{\{\mathcal{M}\}}^{\on{b}}$. To see this we show the following result; for the case $\mathcal{M}:=\{M\}$ see \cite[Example 4.4.]{KMRplot}.

\begin{lemma}\label{counterexampleunionspace2}
Let $\mathcal{M}=\{M^x: x\in\Lambda\}$ be \hyperlink{Msc}{$(\mathcal{M}_{\on{sc}})$} with $\Lambda=\NN_{>0}$.

Then there exists $f:\RR^2\rightarrow\RR^{\NN_{>0}}$ which is $\mathcal{E}_{\{\mathcal{M}\}}$, but there is no reasonable topology on $\mathcal{E}_{\{\mathcal{M}\}}(\RR,\RR^{\NN_{>0}})$ such that the associated mapping $f^{\vee}:\RR\rightarrow\mathcal{E}_{\{\mathcal{M}\}}(\RR,\RR^{\NN_{>0}})$ is $\mathcal{E}_{\{\mathcal{M}\}}^{\on{b}}$.
\end{lemma}
For a ''reasonable topology'' on $\mathcal{E}_{\{\mathcal{M}\}}(\RR,\RR^{\NN_{>0}})$ we assume only that all point-evaluations $\ev_t:\mathcal{E}_{\{\mathcal{M}\}}(\RR,\RR^{\NN_{>0}})\rightarrow\RR^{\NN_{>0}}$ are bounded linear mappings.

\demo{Proof}
Consider $f:\RR^2\rightarrow\RR^{\NN_{>0}}$ defined by $f(s,t):=(\theta_x(st))_{x\in\Lambda}$, $\theta_x\in\mathcal{E}^{\on{global}}_{\{M^x\}}(\RR,\RR)$, see \hyperlink{chf}{$(\text{chf})$}. $f$ is clearly $\mathcal{E}_{\{\mathcal{M}\}}$ since each linear functional on $\RR^{\NN_{>0}}$ depends only on finitely many coordinates. If $f^{\vee}:\RR\rightarrow\mathcal{E}_{\{\mathcal{M}\}}(\RR,\RR^{\NN_{>0}})$ would be $\mathcal{E}_{\{\mathcal{M}\}}^{\on{b}}$, then there would exist $h>0$ and some $y\in\Lambda$ such that the set
$$\left\{\frac{(f^{\vee})^{(k)}(0)}{h^k M^y_k}: k\in\NN\right\}$$
would be bounded in $\mathcal{E}_{\{\mathcal{M}\}}(\RR,\RR^{\NN_{>0}})$. But if we apply the bounded linear function $\ev_t$ for $t=2h$, then
$$\frac{|(f^{\vee})^{(k)}(0)(2h)|}{h^k M^y_k}=\left(\frac{(2h)^k|\theta_x^{(k)}(0)|}{h^k M^y_k}\right)_{x\in\Lambda}\ge\left(\frac{2^k M^x_k}{M^y_k}\right)_{x\in\Lambda}$$
and so the coordinates are unbounded as $k\rightarrow\infty$ whenever $x\ge y$.
\qed\enddemo

To get $\mathcal{E}_{\{\mathcal{M}\}}=\mathcal{E}_{\{\mathcal{M}\}}^{\on{b}}$ we have to assume additional assumptions, see \cite[Lemma 4.3.]{KMRplot} for the constant case.

\begin{lemma}\label{Plots-existencequantifiers}
Let $\mathcal{M}$ be \hyperlink{Marb}{$(\mathcal{M})$}, let $E,F$ be convenient and let $U\subseteq E$ be a $c^{\infty}$-open subset. Assume that there exists a Baire-vector-space-topology on the dual $F^{\ast}$ for which the point evaluations $\ev_x$ are continuous for all $x\in F$. Then $f:U\rightarrow F$ is $\mathcal{E}_{\{\mathcal{M}\}}$ if and only if $f$ is $\mathcal{E}_{\{\mathcal{M}\}}^{\on{b}}$.
\end{lemma}

\demo{Proof}
$(\Leftarrow)$ is clear.

$(\Rightarrow)$ Let $B$ a closed absolutely convex bounded subset of $E$, furthermore consider a compact set $K$ in $U_B$ (w.r.t. $\|\cdot\|_B$) and introduce the sets
$$A_{x,h,C}:=\left\{\alpha\in F^{\ast}:\frac{|(\alpha\circ f)^{(k)}(a)(v_1,\dots,v_k)|}{h^k M^x_k}\le C,\;\forall\;k\in\NN, a\in K, \|v_i\|_B\le 1\right\}.$$
These sets are closed in $F^{\ast}$ for the Baire-topology and $\bigcup_{x\in\Lambda,h,C>0}A_{x,h,C}=F^{\ast}$ holds. Then, by the {\itshape Baire-property} of $F^{\ast}$, there exist $x_0\in\Lambda$, $h_0,C_0>0$ such that the interior $\overset{\circ}{A}_{x_0,h_0,C_0}$ is non-empty. Let $\alpha_0\in \overset{\circ}{A}_{x_0,h_0,C_0}$, then for all $\alpha\in F^{\ast}$ there exists $\varepsilon>0$, such that we get $\varepsilon\alpha\in \overset{\circ}{A}_{x_0,h_0,C_0}-\alpha_0\Leftrightarrow\varepsilon\alpha+\alpha_0\in\overset{\circ}{A}_{x_0,h_0,C_0}$.

Thus for all $a\in K$, $k\in\NN$ and $\|v_i\|_B\le 1$ we get
\begin{align*}
|(\alpha\circ f)^{(k)}(a)(v_1,\dots,v_k)|&\le\frac{1}{\varepsilon}\left(|((\varepsilon\alpha)+\alpha_0)\circ f)^{(k)}(a)|+|(\alpha_0\circ f)^{(k)}(a)|\right)\le\frac{2 C_0}{\varepsilon} h_0^k M^{x_0}_k.
\end{align*}
So the set $\left\{\frac{f^{(k)}(a)(v_1,\dots,v_k)}{h_0^k M^{x_0}_k}: k\in\NN, a\in K, \|v_i\|_B\le 1\right\}$ is weakly bounded (in $F$), hence bounded. Since $B$ was arbitrary we get $f\in\mathcal{E}_{\{\mathcal{M}\}}^{\on{b}}$.
\qed\enddemo

If the matrix is non-constant and has infinite index set, e.g. if $\mathcal{M}$ is coming from $\omega\in\hyperlink{omset1}{\mathcal{W}}$ which does not have \hyperlink{om6}{$(\omega_6)$} - see \cite[Section 5]{compositionpaper}, then another phenomenon appears.

\begin{proposition}\label{counterexampleunionspace}
Let $\mathcal{M}=\{M^x: x\in\Lambda=\NN_{>0}\}$ be \hyperlink{Msc}{$(\mathcal{M}_{\on{sc}})$} with \hyperlink{R-strict}{$(\mathcal{M}_{\{\on{strict}\}})$}.

Then there exist locally convex vector spaces $E$ and $\mathcal{E}_{\{\mathcal{M}\}}$-curves $c:\RR\rightarrow E$ that are not $\mathcal{E}_{\{M^x\}}$ for any $x\in\Lambda$, i.e. $\mathcal{E}_{\{\mathcal{M}\}}(\RR,E)\subsetneq\bigcup_{x\in\Lambda}\mathcal{E}_{\{M^x\}}(\RR,E)$.
\end{proposition}

\demo{Proof}
By \hyperlink{R-strict}{$(\mathcal{M}_{\{\text{strict}\}})$} we have that for each $x\in\Lambda$ we can find $x_1\in\Lambda$, $x_1>x$, such that $\mathcal{E}_{\{M^{x}\}}\subsetneq\mathcal{E}_{\{M^{x_1}\}}$. Iterating \hyperlink{R-strict}{$(\mathcal{M}_{\{\on{strict}\}})$} we obtain a strictly increasing sequence $(x_i)_{i\ge 0}$ with $x_0=x$ and $\lim_{i\rightarrow\infty}x_i=+\infty$, w.l.o.g. one could assume that $\mathcal{M}=\{M^{x_i}: i\in\NN\}$.

So let $x\in\Lambda$ be arbitrary but from now on fixed and set $E:=\RR^{\NN}$. Consider a curve $c:\RR\rightarrow\RR^{\NN}$, $c(t)=(c_i(t))_{i\in\NN}=(c_0(t),c_1(t),\dots)$, with the following property: $c_0$ is $\mathcal{E}_{\{M^{x_0}\}}^{\on{b}}$, and for each $i\ge 1$ we assume $c_i\in\mathcal{E}_{\{M^{x_i}\}}\backslash\mathcal{E}_{\{M^{x_{i-1}}\}}$.

The curve $c$ is $\mathcal{E}_{\{\mathcal{M}\}}$ since each $\alpha\in(\RR^{\NN})^{\ast}=\RR^{(\NN)}$ depends only on finitely many coordinates. Let $i$ be the maximal of these coordinates. Then $\alpha\circ c\in\mathcal{E}_{\{M^{x_i}\}}(\RR,\RR)$, thus $c\in\mathcal{E}_{\{\mathcal{M}\}}(\RR,\RR^{\NN})$.

If there would exist some $y\in\Lambda$ such that $c$ is $\mathcal{E}_{\{M^y\}}$, then for each $\alpha\in\RR^{(\NN)}$ we would get that $\alpha\circ c\in\mathcal{E}_{\{M^y\}}(\RR,\RR)$. According to this $y$ we choose a linear functional $\alpha$ depending on at least $i_0+1$ many coordinates where $x_{i_0}>y$.
\qed\enddemo

\subsection{Closedness under composition of $\mathcal{E}_{[\mathcal{M}]}$}\label{subsection42}

\begin{definition}
Let $E$ be a convenient vector space. A {\itshape $\mathcal{E}_{[\mathcal{M}]}$-Banach-plot} in $E$ is a mapping $c:D\rightarrow E$ such that $c\in\mathcal{E}_{[\mathcal{M}]}$ and $D$ denotes an open set in some Banach space $F$. It is sufficient to consider the open unit ball $D=oF$.
\end{definition}

Using the definitions and projective representations of section \ref{section3} we can generalize \cite[Theorem 4.8.]{KMRplot}.

\begin{theorem}\label{Banach-plot-theorem}
Let $\mathcal{M}$ be \hyperlink{Msc}{$(\mathcal{M}_{\on{sc}})$} with $\Lambda=\RR_{>0}$, let $U\subseteq E$ be a $c^{\infty}$-open subset in a convenient vector space $E$ and $F$ be a Banach space.

If $\mathcal{M}$ has $(\mathcal{M}_{[\on{FdB}]})$ and $f:U\rightarrow F$, then $f\in\mathcal{E}_{[\mathcal{M}]}$ implies $f\circ c\in\mathcal{E}_{[\mathcal{M}]}$ for all $\mathcal{E}_{[\mathcal{M}]}$-Banach plots $c$.

The converse implication holds always by the definitions given in \ref{subsection22}.
\end{theorem}

\demo{Proof}
We follow the proof of \cite[Theorem 4.8.]{KMRplot} and apply Proposition \ref{existencequantorroumieustrong} for the Roumieu- and Proposition \ref{existencequantorbeurlingstrong} for the Beurling-case.

$(a)$ Beurling-case $\mathcal{E}_{(\mathcal{M})}$.

We have to show that $f\circ c$ is $\mathcal{E}_{(\mathcal{M})}$ for each $\mathcal{E}_{(\mathcal{M})}$-Banach-plot $c:G\supseteq D\rightarrow E$, where $D$ denotes the open unit ball in an arbitrary Banach-space $G$. By $(3)$ in Proposition \ref{existencequantorbeurlingstrong} we have to prove that for each compact $K\subseteq D$ and for each $(r_k)_k\in\hyperlink{r-beur-sub}{\mathcal{R}_{\text{Beur},\text{sub}}}$, $(s_k)_k\in\hyperlink{s-beur-FdB}{\mathcal{S}^{\mathcal{M}}_{\text{Beur},\text{FdB}}}$ the set
$$\left\{\frac{(f\circ c)^{(k)}(a)(v_1,\dots,v_k)}{k!} r_k s_k\delta^k: a\in K, k\in\NN, \|v_i\|_E\le 1\right\}$$
is bounded in $F$ for each $\delta>0$. So let $\delta>0$, the sequences $(r_k)_k$, $(s_k)_k$, and finally a compact (w.l.o.g. convex) set $K\subseteq D$ be given, arbitrary but from now on fixed. Then for each $\alpha\in E^{\ast}$ by assumption and by $(2)$ in Proposition \ref{existencequantorbeurlingstrong} applied to the sequence $(r_k(2 D\delta)^k)_k$ and $(\hat{s}_k)_k\in\hyperlink{S-beur}{\mathcal{S}^{\mathcal{M}}_{\on{Beur}}}$, where the constant $D$ is coming from $s_k\le D^k(\hat{s}_o)_k$ (since $(s_k)_k\in\hyperlink{s-beur-FdB}{\mathcal{S}^{\mathcal{M}}_{\text{Beur},\text{FdB}}}$), the set
\begin{equation}\label{plots-set}
\left\{\frac{(\alpha\circ c)^{(k)}(a)(v_1,\dots,v_k) r_k\hat{s}_k(2 D\delta)^k}{k!}: a\in K, k\in\NN, \|v_i\|_G\le 1\right\}
\end{equation}
is bounded in $\RR$. So the set
$$\left\{\frac{c^{(k)}(a)(v_1,\dots,v_k)r_k\hat{s}_k(2 D\delta)^k}{k!}: a\in K, k\in\NN, \|v_i\|_G\le 1\right\}$$
is contained in some closed absolutely convex bounded subset $B$ of $E$, hence
\begin{equation}\label{plots-set1}
\frac{\|c^{(k)}(a)\|_{L^k(G,E_B)} r_k\hat{s}_k\delta^k}{k!}\le\frac{1}{(2D)^k}.
\end{equation}
We proceed now as in \cite[Theorem 4.8.]{KMRplot}. $c(K)$ is compact in $E_B$ since the mapping $c:K\rightarrow E_B$ is Lipschitzian: For all $a,b\in K$ we get $c(a)-c(b)\in\frac{\|a-b\|_G}{2Dr_1\hat{s}_1\delta}B$. Then we estimate for all $\delta>0$ and $k\in\NN_{>0}$ as follows:
\begin{align*}
&\left\|\frac{(f\circ c)^{(k)}(a)}{k!}r_ks_k\delta^k\right\|_{L^k(G,F)}
\\&
\le\sum_{j\ge 0}\sum_{\alpha\in\NN^j_{>0},\sum_{i=1}^j\alpha_i=k}D^k\frac{\left\| f^{(j)}(c(a))\right\|_{L^j(E_B,F)}\hat{s}_j}{j!}\underbrace{\prod_{i=1}^j\frac{\|c^{(\alpha_i)}(a)\|_{L^{\alpha_i}(G,E_B)} r_{\alpha_i}\hat{s}_{\alpha_i}\delta^{\alpha_i}}{\alpha_i!}}_{\le\frac{1}{(2D)^{\alpha_1}}\cdots\frac{1}{(2D)^{\alpha_j}}=\frac{1}{(2D)^k}}
\\&
\le\left(\frac{1}{2}\right)^k\sum_{j\ge 0}\sum_{\alpha\in\NN^j_{>0},\sum_{i=1}^j\alpha_i=k}\underbrace{\frac{\left\| f^{(j)}(c(a))\right\|_{L^j(E_B,F)}}{j! m^x_j}}_{(\star)\le C h^j}\underbrace{(\hat{s}_j m^x_j)}_{\le C_1^j}
\\&
\le(C h C_1)\left(\frac{1}{2}\right)^k\cdot\sum_{j\ge 0}\binom{k-1}{j-1}(h C_1)^{j-1}=(C h C_1)\left(\frac{1}{2}\right)^k(1+C_1 h)^{k-1}
\\&
\le(C h C_1)\left(\frac{(1+C_1 h)}{2}\right)^k.
\end{align*}
We have to choose $x\in\Lambda$ according to $(\hat{s}_j)_j\in\hyperlink{s-beur}{\mathcal{S}^{\mathcal{M}}_{\text{Beur}}}$ (arising in \hyperlink{s-beur-FdB}{$\mathcal{S}^{\mathcal{M}}_{\text{Beur},\text{FdB}}$}) such that $\hat{s}_j m^x_j\le C_1^j$ for some constant $C_1>0$ and all $j\in\NN$. Since $f\in\mathcal{E}_{(\mathcal{M})}$, we obtain the estimate $(\star)$ with this index $x$ and arbitrary $h>0$ for a constant $C=C_{x,h}$ and all $j\in\NN$. Finally we can choose $h:=\frac{1}{C_1}$ and so the expression at the beginning is bounded by $C=C_{x,1/C_1}$.\vspace{6pt}

$(b)$ Roumieu-case $\mathcal{E}_{\{\mathcal{M}\}}$.

Use Proposition \ref{existencequantorroumieustrong} and by $(3)$ there it is sufficient to show that each compact $K\subseteq D$ and for each $(r_k)_k\in\hyperlink{r-roum-sub}{\mathcal{R}_{\text{Roum},\text{sub}}}$, $(s_k)_k\in\hyperlink{s-roum-FdB}{\mathcal{S}^{\mathcal{M}}_{\text{Roum},\text{FdB}}}$, there exists $\varepsilon>0$ such that the set
\begin{equation}\label{plots-set2}
\left\{\frac{(f\circ c)^{(k)}(a)(v_1,\dots,v_k)}{k!} r_k s_k\varepsilon^k: a\in K, k\in\NN, \|v_i\|_E\le 1\right\}
\end{equation}
is bounded in $F$.

We use the same proof as above and replace in $(2)$ in Proposition \ref{existencequantorroumieustrong} the sequence $(r_k)_k$ by $((2D)^k r_k)_k$, where $D$ is the constant arising in $s_k\le D^k(\hat{s}_o)_k$ (since $(s_k)_k\in\hyperlink{s-roum-FdB}{\mathcal{S}^{\mathcal{M}}_{\text{Roum},\text{FdB}}}$ and so $(\hat{s}_k)_k\in\hyperlink{S-roum}{\mathcal{S}^{\mathcal{M}}_{\on{Roum}}}$). Then we take $\delta=1$ in \eqref{plots-set}, in \eqref{plots-set1} and in the Lipschitz-argument. We can use now precisely the same estimate as for the Beurling-case (for $\delta=1$) and so we have shown \eqref{plots-set2} for $\varepsilon=\frac{2}{(1+C_1 h)}$. Note that $f\in\mathcal{E}_{\{\mathcal{M}\}}$, hence we have to consider $x\in\Lambda$ and $h>0$ sufficiently large to obtain estimate $(\star)$ for some constant $C$. According to this chosen $x\in\Lambda$ we can estimate $\hat{s}_j m^x_j\le C_1^j$ for a constant $C_1$ and all $j\in\NN$, since $(\hat{s}_j)_j\in\hyperlink{s-roum}{\mathcal{S}^{\mathcal{M}}_{\text{Roum}}}$.
\qed\enddemo

Using Theorem \ref{Banach-plot-theorem} we can generalize \cite[Theorem 4.9.]{KMRplot}.

\begin{theorem}\label{Banach-plots-category}
Let $\mathcal{M}$ be \hyperlink{Msc}{$(\mathcal{M}_{\on{sc}})$} with $\Lambda=\RR_{>0}$. Let $E,F,G$ be convenient vector spaces, $U\subseteq E$ and $V\subseteq F$ be $c^{\infty}$-open and $f:U\rightarrow F$, $g: V\rightarrow G$ with $f(U)\subseteq V$.

\begin{itemize}
\item[$(a)$] If \hyperlink{B-FdB}{$\mathcal{M}_{(\on{FdB})}$}, then $f,g\in\mathcal{E}_{(\mathcal{M})}$ implies $g\circ f\in\mathcal{E}_{(\mathcal{M})}$.

\item[$(b)$] If \hyperlink{R-FdB}{$\mathcal{M}_{\{\on{FdB}\}}$}, then $f,g\in\mathcal{E}_{\{\mathcal{M}\}}$ implies $g\circ f\in\mathcal{E}_{\{\mathcal{M}\}}$.
\end{itemize}
\end{theorem}

\demo{Proof}
By definition of $\mathcal{E}_{[\mathcal{M}]}$ we have to show that for all closed absolutely convex bounded subsets $B\subseteq E$ and for all $\alpha\in G^{\ast}$ the composite $\alpha\circ g\circ f\circ i_B: U_B\rightarrow\RR$ is $\mathcal{E}_{[\mathcal{M}]}$. By assumption $f\circ i_B\in\mathcal{E}_{[\mathcal{M}]}$ and $\alpha\circ g\in\mathcal{E}_{[\mathcal{M}]}$ hold, so we can use Theorem \ref{Banach-plot-theorem} to obtain the desired implication. Note that $f\circ i_B$ is a $\mathcal{E}_{[\mathcal{M}]}$-Banach plot.
\qed\enddemo

\section{Exponential laws for $\mathcal{E}_{[\mathcal{M}]}$}\label{section5}

We start with the generalization of \cite[Lemma 5.1.]{KMRplot}.

\begin{lemma}\label{generalexplemma}
Let $\mathcal{M}$ be \hyperlink{Marb}{$(\mathcal{M})$} or \hyperlink{Msc}{$(\mathcal{M}_{\on{sc}})$} with $\Lambda=\RR_{>0}$, let $E$ be Banach and $U\subseteq E$ open. Let $F$ be convenient and $\mathcal{B}$ a family of bounded linear functionals on $F$ which together detect bounded sets, i.e. $B\subseteq E$ is bounded in $E$ if and only if $\alpha(B)$ is bounded in $\RR$ for all $\alpha\in\mathcal{B}$. Then we have
$$f\in\mathcal{E}_{[\mathcal{M}]}(U,F)\Leftrightarrow\alpha\circ f\in\mathcal{E}_{[\mathcal{M}]}(U,\RR)\;\;\;\hspace{20pt}\;\forall\;\alpha\in\mathcal{B}.$$
\end{lemma}
\demo{Proof}
For $\mathcal{E}$-curves this follows by \cite[2.1., 2.11.]{KM97}, and so by composing with such curves for $\mathcal{E}$-mappings $f:U\rightarrow F$.

In the Roumieu-case we use $(1)\Leftrightarrow(2)$ in Proposition \ref{existencequantorroumieuweak}. Hence for arbitrary $\alpha\in F^{\ast}$ the mapping $\alpha\circ f$ is $\mathcal{E}_{\{\mathcal{M}\}}$ if and only if for each compact $K\subseteq U$ the set
$$\left\{(\alpha\circ f)^{(k)}(a)(v_1,\dots,v_k) r_k s_k: a\in K, k\in\NN, \|v_i\|_E\le 1\right\}$$
is bounded in $\RR$ for each $(r_k)_k\in\hyperlink{r-roum}{\mathcal{R}_{\text{Roum}}}$ and for each $(s_k)_k\in\hyperlink{S-roum}{\tilde{\mathcal{S}}^{\mathcal{M}}_{\text{Roum}}}$. So the smooth mapping
$f:U\rightarrow F$ is $\mathcal{E}_{\{\mathcal{M}\}}$ if and only if the set
$$\left\{f^{(k)}(a)(v_1,\dots,v_k) r_k s_k: a\in K, k\in\NN, \|v_i\|_E\le 1\right\}$$
is bounded in $F$, for each compact $K\subseteq U$, $(r_k)_k\in\hyperlink{r-roum}{\mathcal{R}_{\text{Roum}}}$ and for each $(s_k)_k\in\hyperlink{S-roum}{\tilde{\mathcal{S}}^{\mathcal{M}}_{\text{Roum}}}$.

Because $\mathcal{B}$ detects bounded sets we can replace in the above equivalences $F^{\ast}$ by $\mathcal{B}$.

For the Beurling-case proceed analogously and use $(1)\Leftrightarrow(2)$ in Proposition \ref{existencequantorbeurlingweak}.
\qed\enddemo

Now we are able to prove {\itshape Cartesian closedness} for classes $\mathcal{E}_{[\mathcal{M}]}$ and so generalize \cite[Theorem 5.2.]{KMRplot}.

\begin{theorem}\label{exponentiallaw}
Let $\mathcal{M}$ be \hyperlink{Msc}{$(\mathcal{M}_{\on{sc}})$} with $\Lambda=\RR_{>0}$, let $U_i\subseteq E_i$ be $c^{\infty}$-open subsets in convenient vector spaces $E_i$ for $i=1,2$ and moreover let $F$ be also a convenient vector space. Then we obtain:

\begin{itemize}
\item[$(a)$] If \hyperlink{R-mg}{$(\mathcal{M}_{\{\on{mg}\}})$}, then
$$f\in\mathcal{E}_{\{\mathcal{M}\}}(U_1\times U_2,F)\Longleftrightarrow f^{\vee}\in\mathcal{E}_{\{\mathcal{M}\}}(U_1,\mathcal{E}_{\{\mathcal{M}\}}(U_2,F)).$$

\item[$(b)$] If \hyperlink{B-mg}{$(\mathcal{M}_{(\on{mg})})$}, then
$$f\in\mathcal{E}_{(\mathcal{M})}(U_1\times U_2,F)\Longleftrightarrow f^{\vee}\in\mathcal{E}_{(\mathcal{M})}(U_1,\mathcal{E}_{(\mathcal{M})}(U_2,F)).$$
\end{itemize}
\end{theorem}
Important remarks:
\begin{itemize}
\item[$(i)$] In both cases $(\Longleftarrow)$ holds also without \hyperlink{R-mg}{$(\mathcal{M}_{\{\on{mg}\}})$} respectively  \hyperlink{B-mg}{$(\mathcal{M}_{(\on{mg})})$}.

\item[$(ii)$] To prove $(\Longleftarrow)$ it is sufficient to assume that $\mathcal{M}$ is \hyperlink{Marb}{$(\mathcal{M})$} and $(\mathcal{M}_{[\on{alg}]})$.

\item[$(iii)$] For the proof it is not necessary to assume that $\mathcal{E}_{\{\mathcal{M}\}}$ respectively $\mathcal{E}_{(\mathcal{M})}$ is a category, i.e. closedness under composition.

\item[$(iv)$] If $\mathcal{M}$ is \hyperlink{Msc}{$(\mathcal{M}_{\on{sc}})$} with $\Lambda=\RR_{>0}$, $(\mathcal{M}_{[\on{mg}]})$ and $(\mathcal{M}_{[\on{FdB}]})$, then by Theorem \ref{exponentiallaw} and Theorem \ref{Banach-plots-category} the category $\mathcal{E}_{[\mathcal{M}]}$ is cartesian closed.
\end{itemize}

\demo{Proof}
The technique and methods are completely analogous to \cite[Theorem 5.2.]{KMRplot}, for convenience of the reader we give the full proof.

As shown in \cite[3.12.]{KM97} we have $\mathcal{E}(U_1\times U_2,F)\cong\mathcal{E}(U_1,\mathcal{E}(U_2,F))$. So we assume form now on that all occurring mappings are smooth. Let $B\subseteq E_1\times E_2$ and $B_i\subseteq E_i$, $i=1,2$, where $B,B_1,B_2$ run through all closed absolutely convex bounded subsets. Similarly as shown in \cite[Theorem 5.2.]{KMRplot} we get:
\begin{align*}
&f\in\mathcal{E}_{[\mathcal{M}]}(U_1\times U_2,F)
\\&
\Leftrightarrow\;\forall\;\alpha\in F^{\ast}\;\forall\;B:\;\alpha\circ f\circ i_B\in\mathcal{E}_{[\mathcal{M}]}((U_1\times U_2)_B,\RR)
\\&
\Leftrightarrow\;\forall\;\alpha\in F^{\ast}\;\forall\;B_1,B_2:\;\alpha\circ f\circ(i_{B_1}\times i_{B_2})\in\mathcal{E}_{[\mathcal{M}]}((U_1)_{B_1}\times(U_2)_{B_2},\RR)
\end{align*}
and
\begin{align*}
&f^{\vee}\in\mathcal{E}_{[\mathcal{M}]}(U_1,\mathcal{E}_{[\mathcal{M}]}(U_2,F))
\\&
\Leftrightarrow\;\forall\;B_1:\;f^{\vee}\circ i_{B_1}\in\mathcal{E}_{[\mathcal{M}]}((U_1)_{B_1},\mathcal{E}_{[\mathcal{M}]}(U_2,F))
\\&
\Leftrightarrow\;\forall\;\alpha\in F^{\ast}\;\forall\;B_1,B_2:\;\mathcal{E}_{[\mathcal{M}]}(i_{B_2},\alpha)\circ f^{\vee}\circ i_{B_1}\in\mathcal{E}_{[\mathcal{M}]}((U_1)_{B_1},\mathcal{E}_{[\mathcal{M}]}((U_2)_{B_2},\RR)),
\end{align*}
where Lemma \ref{generalexplemma} is used and note that the linear mappings $\mathcal{E}_{[\mathcal{M}]}(i_{B_2},\alpha)$ generate the bornology.

With these preparations we are able to restrict ourselves to $U_i\subseteq E_i$ open sets in Banach spaces $E_i$ and $F=\RR$. We start now with $(\Longrightarrow)$ for both cases.\vspace{6pt}

Let $f\in\mathcal{E}_{[\mathcal{M}]}(U_1\times U_2,\RR)$, then clearly $f^{\vee}$ takes values in the space $\mathcal{E}_{[\mathcal{M}]}(U_2,\RR)$.

First we show that

{\itshape Claim.} $f^{\vee}: U_1\rightarrow\mathcal{E}_{[\mathcal{M}]}(U_2,\RR)$ is $\mathcal{E}$ with $d^jf^{\vee}=(\partial^j_1 f)^{\vee}$.

$\mathcal{E}_{[\mathcal{M}]}(U_2,\RR)$ are convenient vector spaces, hence by \cite[5.20.]{KM97} it suffices to prove that the iterated unidirectional derivatives $d^j_vf^{\vee}(x)$ exist, are equal to $\partial^j_1 f(x,\cdot)(v^j)$, and are separately bounded for $x$ and $v$ in compact subsets. For $j=1$ and $x,v,y$ fixed we consider the smooth curve $c:t\mapsto f(x+tv,y)$. Then, by the fundamental theorem of calculus, we obtain:
\begin{align*}
&\frac{f^{\vee}(x+tv)-f^{\vee}(x)}{t}(y)-(\partial_1 f)^{\vee}(x)(y)(v)=\frac{c(t)-c(0)}{t}-c'(0)
\\&
=t\int_0^1s\int_0^1c''(tsr)drds=t\int_0^1s\int_0^1\partial^2_1f(x+tsrv,y)(v,v)drds.
\end{align*}
$(\partial^2_1f)^{\vee}(K_1)(o(E_1\times E_1))$ is bounded in $\mathcal{E}_{[\mathcal{M}]}(U_2,\RR)$ and for each compact set $K_1\subseteq U_1$ this expression is Mackey-convergent to $0$ in $\mathcal{E}_{[\mathcal{M}]}(U_2,\RR)$ as $t\rightarrow 0$. Hence $d_vf^{\vee}(x)$ exists an is equal to $\partial_1 f(x,\cdot)(v)$.

The induction argument is completely the same as in \cite[Theorem 5.2.]{KMRplot}.\vspace{6pt}

We distinguish now between the Roumieu- and the Beurling-case.

{\itshape The Beurling-case.}

We have to show that $f^{\vee}:U_1\rightarrow\mathcal{E}_{(\mathcal{M})}(U_2,\RR)$ is $\mathcal{E}_{(\mathcal{M})}$.

By Lemma \ref{generalexplemma} it suffices to prove that $f^{\vee}:U_1\rightarrow\mathcal{E}_{M^x,h}(E_2\supseteq K_2,\RR)$ is $\mathcal{E}_{(\mathcal{M})}^{\on{b}}=\mathcal{E}_{(\mathcal{M})}$ for each $K_2\subseteq U_2$ compact, each $h>0$ and $x\in\Lambda=\RR_{>0}$. This holds, because each $\alpha\in(\mathcal{E}_{(\mathcal{M})}(U_2,\RR))^{\ast}$ factorizes over $\mathcal{E}_{M^x,h}(E_2\supseteq K_2,\RR)$ for some $K_2$, $h$ and $x$.

So we have to show that for each compact sets $K_1\subseteq U_1$, $K_2\subseteq U_2$, each $h_1,h_2>0$ and each $x_1,x_2\in\Lambda$, the set
\begin{equation}\label{exponentiallawequ1}
\left\{\frac{d^{k_1}f^{\vee}(a_1)(v_1^1,\dots,v^1_{k_1})}{h_1^{k_1} M^{x_1}_{k_1}}: a_1\in K_1, k_1\in\NN, \|v^1_j\|_{E_1}\le 1\right\}
\end{equation}
is bounded in the space $\mathcal{E}_{M^{x_2},h_2}(E_2\supseteq K_2,\RR)$. Equivalently, for all compact sets $K_1,K_2$, for all $h_1,h_1>0$ and all $x_1,x_2\in\Lambda$ the set
\begin{equation}\label{exponentiallawequ2}
\left\{\frac{\partial^{k_2}_2\partial^{k_1}_1 f(a_1,a_2)(v^1_1,\dots,v^1_{k_1};v^2_1,\dots,v^2_{k_2})}{h_2^{k_2} h_1^{k_1} M^{x_2}_{k_2} M^{x_1}_{k_1}}: a_i\in K_i, k_i\in\NN, \|v^i_j\|_{E_i}\le 1; i=1,2\right\}
\end{equation}
is bounded in $\RR$.

Let $a_1\in K_1$, $k_1\in\NN$, then we obtain the following estimate:
\begin{align*}
&\left\|\frac{d^{k_1}f^{\vee}(a_1)(v_1^1,\dots,v^1_{k_1})}{h_1^{k_1} M^{x_1}_{k_1}}\right\|^J_{M^{x_2},K_2,h_2}
\\&
=\sup\left\{\frac{\left|\partial^{k_2}_2\partial^{k_1}_1f(a_1,a_2)(v^1_1,\dots,v^1_{k_1};v^2_1,\dots,v^2_{k_2})\right|}{h_1^{k_1} h_2^{k_2} M^{x_1}_{k_1} M^{x_2}_{k_2}}: a_2\in K_2, k_2\in\NN, \|v^2_j\|_{E_2}\le 1\right\}
\\&
\underbrace{\le}_{\hyperlink{B-mg}{(\mathcal{M}_{(\text{mg})})}}\sup\left\{C^{k_1+k_2}\frac{\left|\partial^{k_2}_2\partial^{k_1}_1f(a_1,a_2)(v^1_1,\dots,v^1_{k_1};v^2_1,\dots,v^2_{k_2})\right|}{h_1^{k_1} h_2^{k_2} M^y_{k_1+k_2}}: a_2\in K_2, k_2\in\NN, \|v^2_j\|_{E_2}\le 1\right\}
\\&
\le\sup\left\{\frac{\left|\partial^{k_2}_2\partial^{k_1}_1f(a_1,a_2)(v^1_1,\dots,v^1_{k_1};v^2_1,\dots,v^2_{k_2})\right|}{h^{k_1+k_2} M^y_{k_1+k_2}}: a_2\in K_2, k_2\in\NN, \|v^2_j\|_{E_2}\le 1\right\}<+\infty,
\end{align*}
where we have put $h:=\frac{1}{C}\min\{h_1,h_2\}$. Note that $f$ is $\mathcal{E}_{(\mathcal{M})}$ and so for arbitrary $h_1,h_2>0$ and $x_1,x_2\in\Lambda$ we can find $y\in\Lambda$ and $h>0$ such that the last inequality is valid. This shows that $f^{\vee}$ is $\mathcal{E}_{(\mathcal{M})}$.\vspace{6pt}

{\itshape The Roumieu-case.}

By Lemma \ref{generalexplemma} it suffices to prove that $f^{\vee}:U_1\rightarrow\underset{x_2\in\Lambda}{\varinjlim}\underset{h_2>0}{\varinjlim}\mathcal{E}_{M^{x_2},h_2}(E_2\supseteq K_2,\RR)$ is $\mathcal{E}_{\{\mathcal{M}\}}^{\on{b}}\subseteq\mathcal{E}_{\{\mathcal{M}\}}$ for each compact set $K_2\subseteq U_2$. This holds because each $\alpha\in(\mathcal{E}_{\{\mathcal{M}\}}(U_2,\RR))^{*}$ factorizes over some $\underset{x_2\in\Lambda}{\varinjlim}\underset{h_2>0}{\varinjlim}\mathcal{E}_{M^{x_2},h_2}(E_2\supseteq K_2,\RR)$.

So we have to prove that for all $K_1\subseteq U_1$, $K_2\subseteq U_2$ compact there exist $h_1>0$ and some $x_1\in\Lambda$ such that the set in \eqref{exponentiallawequ1} is bounded in $\underset{x_2\in\Lambda}{\varinjlim}\underset{h_2>0}{\varinjlim}\mathcal{E}_{M^{x_2},h_2}(E_2\supseteq K_2,\RR)$. Equivalently, we have to show that for all $K_1,K_2$ compact there exist $h_1,h_2>0$ and $x_1,x_2\in\Lambda$ such that the set in \eqref{exponentiallawequ2} is bounded in $\RR$.

We can use now the same estimate as for the above Beurling-case and use \hyperlink{R-mg}{$(\mathcal{M}_{\{\text{mg}\}})$}. First, because $f$ is $\mathcal{E}_{\{\mathcal{M}\}}$ and by $(3)$ in Proposition \ref{matrixcompletness} we obtain that there exist some $h>0$ and $y\in\Lambda$, such that the last set

$\sup\left\{\frac{\left|\partial^{k_2}_2\partial^{k_1}_1f(a_1,a_2)(v^1_1,\dots,v^1_{k_1};v^2_1,\dots,v^2_{k_2})\right|}{h^{k_1+k_2} M^y_{k_1+k_2}}: a_2\in K_2, k_2\in\NN, \|v^2_j\|_{E_2}\le 1\right\}$

in the Beurling estimate is bounded. For this $y\in\Lambda$ we obtain by \hyperlink{R-mg}{$(\mathcal{M}_{\{\text{mg}\}})$} that there exist some $x_1,x_2\in\Lambda$ and $C>0$ such that $M^y_{j+k}\le C^{j+k} M^{x_1}_j M^{x_2}_k$ holds for all $j,k\in\NN$. So we can put in the estimate now $h_i:=C h$ for $i=1,2$ to get, that $f^{\vee}$ is $\mathcal{E}_{\{\mathcal{M}\}}$.\vspace{6pt}

Now we start with $(\Longleftarrow)$ for both cases.\vspace{6pt}

Let $f^{\vee}:U_1\rightarrow\mathcal{E}_{[\mathcal{M}]}(U_2,\RR)$ be $\mathcal{E}_{[\mathcal{M}]}$. The mapping $f^{\vee}:U_1\rightarrow\mathcal{E}_{[\mathcal{M}]}(U_2,\RR)\rightarrow\mathcal{E}(U_2,\RR)$ is $\mathcal{E}$, hence it remains to show that $f\in\mathcal{E}_{[\mathcal{M}]}(U_1\times U_2,\RR)$.\vspace{6pt}

{\itshape The Beurling-case.}

For each compact $K_2\subseteq U_2$, each $h_2>0$ and each $x_2\in\Lambda$, the mapping $f^{\vee}:U_1\rightarrow\mathcal{E}_{M^{x_2},h_2}(E_2\supseteq K_2,\RR)$ is $\mathcal{E}_{(\mathcal{M})}^{\on{b}}=\mathcal{E}_{(\mathcal{M})}$. This means that for all compact $K_1\subseteq U_1$, $K_2\subseteq U_2$, each $h_1,h_2>0$ and each $x_1,x_2\in\Lambda$ the set in \eqref{exponentiallawequ1} is bounded in $\mathcal{E}_{M^{x_2},h_2}(E_2\supseteq K_2,\RR)$. Because it is contained in the space $\mathcal{E}_{M^{x_2},K_2,h_2}(U_2,\RR):=\{f\in\mathcal{E}(U_2,\RR): j^{\infty}(f)|_{K_2}\in\mathcal{E}_{M^{x_2},h_2}(E_2\supseteq K_2,\RR)\}$ with semi-norm $\|f\|^J_{M^{x_2},K_2,h_2}:=\|j^{\infty}(f)|_{K_2}\|^J_{M^{x_2},h_2}$, it is also bounded in this space and so the set in \eqref{exponentiallawequ2} is bounded in $\RR$.

By assumption each $M^x$ is log-convex and so $M^x_j M^x_k\le M^x_{j+k}$ for all $j,k\in\NN$. For the next estimate \hyperlink{B-alg}{$(\mathcal{M}_{(\text{alg})})$} would be sufficient. Let $a_1\in K$, $k_1\in\NN$ and $\|v^1_j\|_{E_1}\le 1$, then:
\begin{align*}
+\infty>&\left\|\frac{d^{k_1}f^{\vee}(a_1)(v_1^1,\dots,v_{k_1}^1)}{h_1^{k_1} M^{x_1}_{k_1}}\right\|^J_{M^{x_2},K_2,h_2}
\\&
=\sup\left\{\frac{\left|\partial_2^{k_2}\partial_1^{k_1}f(a_1,a_2)(v_1^1,\dots,v_{k_1}^1;v_1^2,\dots,v^2_{k_2})\right|}{h_1^{k_1} h_2^{k_2} M^{x_1}_{k_1} M^{x_2}_{k_2}}: a_2\in K_2, k_2\in\NN, \left\|v^2_j\right\|_{E_2}\le 1\right\}
\\&
\ge\sup\left\{\frac{\left|\partial_2^{k_2}\partial_1^{k_1}f(a_1,a_2)(v_1^1,\dots,v_{k_1}^1;v_1^2,\dots,v^2_{k_2})\right|}{h^{k_1+k_2} M^y_{k_1+k_2}}: a_2\in K_2, k_2\in\NN, \left\|v^2_j\right\|_{E_2}\le 1\right\}
\end{align*}
where we have put $y:=\max\{x_1,x_2\}$ and $h:=\max\{h_1,h_2\}$ (put $h:=C\max\{h_1,h_2\}$, where $y\in\Lambda$ and $C>0$ are coming from \hyperlink{B-alg}{$(\mathcal{M}_{(\text{alg})})$}). So we have shown that $f$ is $\mathcal{E}_{(\mathcal{M})}$.\vspace{6pt}

{\itshape The Roumieu-case.}

For each compact $K_2\subseteq U_2$ the mapping $f^{\vee}:U_1\rightarrow\underset{x_2\in\Lambda}{\varinjlim}\underset{h_2>0}{\varinjlim}\mathcal{E}_{M^{x_2},h_2}(E_2\supseteq K_2,\RR)$ is $\mathcal{E}_{\{\mathcal{M}\}}$. By $(3)$ in Proposition \ref{matrixcompletness} the dual space $\big(\underset{x_2\in\Lambda}{\varinjlim}\underset{h_2>0}{\varinjlim}\mathcal{E}_{M^{x_2},h_2}(E_2\supseteq K_2,\RR)\big)^{\ast}$ can be equipped with the Baire-vector-space-topology of the countable limit of Banach spaces $\underset{x_2\in\Lambda}{\varprojlim}\underset{h_2>0}{\varprojlim}\left(\mathcal{E}_{M^{x_2},h_2}(E_2\supseteq K_2,\RR)\right)^{\ast}$.

Now we can use Lemma \ref{Plots-existencequantifiers} to conclude that the mapping $f^{\vee}:U_1\rightarrow\underset{x_2\in\Lambda}{\varinjlim}\underset{h_2>0}{\varinjlim}\mathcal{E}_{M^{x_2},h_2}(E_2\supseteq K_2,\RR)$ is $\mathcal{E}_{\{\mathcal{M}\}}^{\on{b}}$.

By $(3)$ in Proposition \ref{matrixcompletness} this inductive limit is countable and compactly regular and so for each compact $K_1\subseteq U_1$ there exist $h_1>0$ and $x_1\in\Lambda$ such that the set in \eqref{exponentiallawequ1} is bounded in $\mathcal{E}_{M^{x_2},h_2}(E_2\supseteq K_2,\RR)$ for some $h_2>0$ and $x_2\in\Lambda$. Because it is contained $\mathcal{E}_{M^{x_2},K_2,h_2}(U_2,\RR):=\{f\in\mathcal{E}(U_2,\RR): j^{\infty}(f)|_{K_2}\in\mathcal{E}_{M^{x_2},h_2}(E_2\supseteq K_2,\RR)\}$ with semi-norm $\|f\|^J_{M^{x_2},K_2,h_2}:=\left\|j^{\infty}(f)|_{K_2}\right\|^J_{M^{x_2},h_2}$, it is also bounded in this space and so the set in \eqref{exponentiallawequ2} is bounded (in $\RR$) with those given $h_1,h_2,x_1,x_2$.

But now we can use the same estimate as in the above Beurling-case to conclude that $f$ is $\mathcal{E}_{\{\mathcal{M}\}}$. Similarly \hyperlink{R-alg}{$(\mathcal{M}_{\{\text{alg}\}})$} would be sufficient for this step.
\qed\enddemo

Using Theorem \ref{exponentiallaw} we can prove now the matrix generalization of \cite[Corollary 5.5.]{KMRplot}:

\begin{corollary}\label{expconsequence1}
Let $\mathcal{M}$ be a weight matrix as assumed in Theorem \ref{exponentiallaw}. Let $E,F,E_i,F_i,G$ be convenient vector spaces and let $U$ and $V$ be $c^{\infty}$-open subsets. Then we get

\begin{itemize}
\item[$(1)$] The exponential law

$\mathcal{E}_{[\mathcal{M}]}(U,\mathcal{E}_{[\mathcal{M}]}(V,G))\cong\mathcal{E}_{[\mathcal{M}]}(U\times V,G)$

holds, it is a linear $\mathcal{E}_{[\mathcal{M}]}$-diffeomorphism of convenient vector spaces.
\end{itemize}

The following mappings are $\mathcal{E}_{[\mathcal{M}]}$:

\begin{itemize}
\item[$(2)$] $\ev: \mathcal{E}_{[\mathcal{M}]}(U,F)\times U\rightarrow F$ given by $\ev(f,x)=f(x)$.

\item[$(3)$] $\ins: E\rightarrow\mathcal{E}_{[\mathcal{M}]}(F,E\times F)$ given by $\ins(x)(y)=(x,y)$.

\item[$(4)$] $(\cdot)^{\wedge}: \mathcal{E}_{[\mathcal{M}]}(U,\mathcal{E}_{[\mathcal{M}]}(V,G))\rightarrow\mathcal{E}_{[\mathcal{M}]}(U\times V,G)$.

\item[$(5)$] $(\cdot)^{\vee}: \mathcal{E}_{[\mathcal{M}]}(U\times V,G)\rightarrow\mathcal{E}_{[\mathcal{M}]}(U,\mathcal{E}_{[\mathcal{M}]}(V,G))$.

\item[$(6)$] $\prod: \prod_i\mathcal{E}_{[\mathcal{M}]}(E_i,F_i)\rightarrow\mathcal{E}_{[\mathcal{M}]}\left(\prod_i E_i,\prod_i F_i\right)$.
\end{itemize}
If $\mathcal{M}$ has also {$(\mathcal{M}_{[\on{FdB}]})$}, then we get
\begin{itemize}
\item[$(7)$] $\comp: \mathcal{E}_{[\mathcal{M}]}(F,G)\times\mathcal{E}_{[\mathcal{M}]}(U,F)\rightarrow\mathcal{E}_{[\mathcal{M}]}(U,G)$.

\item[$(8)$] $\mathcal{E}_{[\mathcal{M}]}(\cdot,\cdot): \mathcal{E}_{[\mathcal{M}]}(F,F_1)\times\mathcal{E}_{[\mathcal{M}]}(E_1,E)\rightarrow\mathcal{E}_{[\mathcal{M}]}(\mathcal{E}_{[\mathcal{M}]}(E,F),\mathcal{E}_{[\mathcal{M}]}(E_1,F_1))$ which is given by $(f,g)\mapsto(h\mapsto f\circ h\circ g)$.
\end{itemize}
\end{corollary}
{\itshape Remark:} $(7)$ proves the claim of \cite[Remark 4.23.]{compositionpaper}.

\subsection{Comparison of conditions $(\text{mg})$ and $(\mathcal{M}_{\{\text{mg}\}})$}
In \cite[Example 5.4.]{KMRplot} it was shown that cartesian closedness fails for $\mathcal{M}=\{M\}$ if $M$ does not satisfy \hyperlink{mg}{$(\on{mg})$}. In the weight matrix case we can prove the following (counter)-example:

\begin{example}\label{Matsumoto-counterexample}
There exist (non-constant) \hyperlink{Msc}{$(\mathcal{M}_{\on{sc}})$} weight matrices $\mathcal{M}$ with \hyperlink{R-mg}{$(\mathcal{M}_{\{\on{mg}\}})$} but such that no $M^x\in\mathcal{M}$ satisfies \hyperlink{mg}{$(\on{mg})$}.
\end{example}

\demo{Proof}
Let $\mathcal{M}=\Omega$ be coming from $\omega\in\hyperlink{omset1}{\mathcal{W}}$ such that \hyperlink{om6}{$(\omega_6)$} does not hold, see \cite[5.5., Corollary 5.8. $(2)$]{compositionpaper}. The weights $\omega(t):=\max\{0,\log(t)^s\}$, $s>1$, are concrete examples, see also \cite{BonetMeiseMelikhov07} for the consequences of \hyperlink{om6}{$(\omega_6)$}.
\qed\enddemo

In the next step we generalize \cite[Example 5.4.]{KMRplot}. We show that \hyperlink{R-mg}{$(\mathcal{M}_{\{\text{mg}\}})$} is necessary for Theorem \ref{exponentiallaw}.

\begin{lemma}\label{moderategrowthcounterexample}
Let $\mathcal{M}$ be \hyperlink{Msc}{$(\mathcal{M}_{\on{sc}})$} with $\Lambda=\NN_{>0}$ but such that \hyperlink{R-mg}{$(\mathcal{M}_{\{\on{mg}\}})$} does not hold. Then there exists $f\in\mathcal{E}_{\{\mathcal{M}\}}(\RR^2,\CC)$ such that the associated mapping $f^{\vee}:\RR\rightarrow\mathcal{E}_{\{\mathcal{M}\}}(\RR,\CC)$ is not $\mathcal{E}_{\{\mathcal{M}\}}$.
\end{lemma}

\demo{Proof}
We follow the proof of \cite[Example 5.4.]{KMRplot}. The negation of \hyperlink{R-mg}{$(\mathcal{M}_{\{\text{mg}\}})$} gives
\begin{equation}\label{R-mg-neg}
\exists\;x\in\Lambda\;\forall\;C>0\;\forall\;y\in\Lambda\;\exists\;j,k\in\NN:\;M^x_{j+k}>C^{j+k} M^y_j M^y_k.
\end{equation}
For this $x\in\Lambda$ and the choice $C=y=n$, $n\in\NN_{>0}$, we obtain sequences $(j_n)_n$ and $(k_n)_n$ such that $(j_n)_n$ is increasing, $j_n\rightarrow\infty$, $k_n\ge 1$ for each $n\in\NN_{>0}$ and with
$$\left(\frac{M^x_{j_n+k_n}}{M^n_{j_n} M^n_{k_n}}\right)^{1/(k_n+j_n)}\ge n.$$
Define a linear functional $\alpha:\mathcal{E}_{\{\mathcal{M}\}}(\RR,\CC)\rightarrow\CC$ by
$$\alpha(f):=\sum_{n\ge 1}(\sqrt{-1})^{3j_n}\frac{f^{(j_n)}(0)}{M^n_{j_n} n^{j_n}}.$$
{\itshape Claim.} $\alpha$ is bounded. For given $f\in\mathcal{E}_{\{\mathcal{M}\}}(\RR,\CC)$ we choose $h>0$ and $l\in\Lambda$ large enough and estimate
$$\left|\sum_{n\ge 0}(\sqrt{-1})^{3j_n}\frac{f^{(j_n)}(0)}{M^n_{j_n} n^{j_n}}\right|\le\sum_{n\ge 0}\frac{\left|f^{(j_n)}(0)\right|}{h^{j_n} M^l_{j_n}}\frac{M^l_{j_n}}{M^n_{j_n}}\left(\frac{h}{n}\right)^{j_n}\le\left\|f\right\|_{M^l,[-1,1],h}\sum_{n\ge 0}\frac{M^l_{j_n}}{M^n_{j_n}}\left(\frac{h}{n}\right)^{j_n}<+\infty.$$
Note that $M^l\le M^n$ for $l\le n$ and $\sum_{n\ge 0}\left(\frac{h}{n}\right)^{j_n}<+\infty$ for each $h>0$.

We apply $\alpha$ to $\tilde{\theta}_x\in\mathcal{E}^{\text{global}}_{\{M^x\}}(\RR,\CC)$ (see \eqref{characcomplex}), where $x\in\Lambda$ is the index from \eqref{R-mg-neg}.

For $s,t\in\RR$ define $\psi_x(s,t):=\tilde{\theta}_x(s+t)$ and so $\psi_x\in\mathcal{E}^{\text{global}}_{\{\mathcal{M}\}}(\RR^2,\CC)$ with $\psi_x^{(\beta_1,\beta_2)}(0,0)=(\sqrt{-1})^{\beta_1+\beta_2} s^x_{\beta_1+\beta_2}$ for all $(\beta_1,\beta_2)\in\NN^2$.

{\itshape Claim.} $\alpha\circ\psi_x^{\vee}$ is not $\mathcal{E}_{\{\mathcal{M}\}}$. Let $h>0$ and $l\in\Lambda$ be arbitrary (large) but fixed and estimate as follows:
\begin{align*}
&\left\|\alpha\circ\psi_x^{\vee}\right\|_{M^l,[-1,1],h}=\sup_{t\in[-1,1],k\in\NN}\frac{|(\alpha\circ\psi_x^{\vee})^{(k)}(t)|}{h^k M^{l}_k}\ge\sup_{k\in\NN}\frac{1}{h^k M^{l}_k}\left|\sum_{n\ge 1}(\sqrt{-1})^{3j_n}\frac{\psi_x^{(j_n,k)}(0,0)}{M^n_{j_n} n^{j_n}}\right|
\\&
=\sup_{k\in\NN}\frac{1}{h^k M^{l}_k}\left|\sum_{n\ge 1}(\sqrt{-1})^{3j_n}\frac{(\sqrt{-1})^{j_n+k} s^x_{j_n+k}}{M^n_{j_n} n^{j_n}}\right|=\sup_{k\in\NN}\frac{1}{h^k M^{l}_k}\left|(\sqrt{-1})^k\sum_{n\ge 1}\frac{s^x_{j_n+k}}{M^n_{j_n} n^{j_n}}\right|
\\&
=\sup_{k\in\NN}\frac{1}{h^k M^{l}_k}\sum_{n\ge 1}\frac{s^x_{j_n+k}}{M^n_{j_n} n^{j_n}}\underbrace{\ge}_{k=k_n}\sup_{n\in\NN_{>0}}\frac{1}{h^{k_n} M^{l}_{k_n}}\frac{M^n_{k_n}}{M^n_{k_n}}\frac{s^x_{j_n+k_n}}{M^n_{j_n} n^{j_n}}
\\&
\ge\sup_{n\in\NN_{>0}}\frac{M^{n}_{k_n}}{h^{k_n} n^{j_n} M^l_{k_n}}\frac{M^x_{j_n+k_n}}{M^n_{j_n} M^{n}_{k_n}}\ge\sup_{n\in\NN_{>0}}\frac{n^{j_n+k_n}}{h^{k_n} n^{j_n}}\frac{M^n_{k_n}}{M^{l}_{k_n}}=+\infty.
\end{align*}
\qed\enddemo

\section{Remarks and special cases}\label{section6}
\subsection{More results for $\mathcal{E}_{[\mathcal{M}]}$}\label{section61}
Let $\mathcal{M}$ be \hyperlink{Marb}{$(\mathcal{M})$} with $\Lambda=\RR_{>0}$. Using the closed graph theorem \cite[52.10]{KM97} the matrix generalization of the uniform boundedness principle \cite[Theorem 6.1.]{KMRplot} is valid for $\mathcal{E}_{[\mathcal{M}]}$, see \cite[Theorem 12.4.1.]{dissertation}. All further results from \cite[Chapter 8]{KMRplot} can be transferred to the matrix-case, see \cite[12.4., 12.6., 12.7.]{dissertation}. For the generalization of \cite[Theorem 2.2.]{KMRplot} see \cite[Proposition 9.4.4.]{dissertation}.\vspace{6pt}

Let $\mathcal{M}$ be \hyperlink{Marb}{$(\mathcal{M})$} and assume that
\begin{itemize}
\item[$(i)$] $\mathcal{M}$ is \hyperlink{Msc}{$(\mathcal{M}_{\text{sc}})$} with $\Lambda=\RR_{>0}$ and has

\item[$(ii)$] $(\mathcal{M}_{[\text{mg}]})$ ($\Rightarrow(\mathcal{M}_{[\text{dc}]})$);

\item[$(iii)$] for the Roumieu-case \hyperlink{holom}{$(\mathcal{M}_{\mathcal{H}})$}, for the Beurling-case \hyperlink{B-Comega}{$(\mathcal{M}_{(\mathcal{C^{\omega}})})$};

\item[$(iv)$] $(\mathcal{M}_{[\text{FdB}]})$ or equivalently $(\mathcal{M}_{[\text{rai}]})$ (see \cite[Lemma 1]{characterizationstabilitypaper}).
\end{itemize}
Using \cite[Theorems 5,6]{characterizationstabilitypaper}, where we characterized the required stability properties for $\mathcal{E}_{[\mathcal{M}]}$, all results from \cite[Chapter 9]{KMRplot} can be transferred to the $\mathcal{E}_{[\mathcal{M}]}$-case, see \cite[12.8.]{dissertation} for full proofs. Note that the characterization theorem for the Beurling-case shown in \cite[Chapter 8]{dissertation} is weaker than \cite[Theorem 6]{characterizationstabilitypaper}.

\subsection{Special cases $\mathcal{M}=\{M\}$ and $\mathcal{M}=\Omega$}\label{section62}
To apply all previous results to the constant case $\mathcal{M}=\{M\}$ we have to assume that
\begin{itemize}
\item[$(i)$] $M\in\hyperlink{LCset}{\mathcal{LC}}$;
\item[$(ii)$] $\liminf_{p\rightarrow\infty}(m_p)^{1/p}>0$ in the Roumieu-, $\lim_{p\rightarrow\infty}(m_p)^{1/p}=+\infty$ in the Beurling-case;
\item[$(iii)$] $M$ has $\hyperlink{mg}{(\text{mg})}(\Rightarrow\hyperlink{dc}{(\text{dc})})$,
\item[$(iv)$] $M$ has \hyperlink{FdB}{$(\text{FdB})$} or equivalently \hyperlink{rai}{$(\text{rai})$} (see also \cite[Chapter 3]{compositionpaper}).
\end{itemize}

If $\mathcal{M}=\Omega=\{(\Omega^l_j)_j: l>0\}$ with $\Omega^l_j:=\exp(1/l\varphi^{*}_{\omega}(lj))$, then we assume that $\omega\in\hyperlink{omset1}{\mathcal{W}}$ and

\begin{itemize}
\item[$(i)$] \hyperlink{om2}{$(\omega_2)$} in the Roumieu-, \hyperlink{om5}{$(\omega_5)$} in the Beurling-case to guarantee \hyperlink{holom}{$(\mathcal{M}_{\mathcal{H}})$} respectively \hyperlink{B-Comega}{$(\mathcal{M}_{(\mathcal{C}^{\omega})})$} (see \cite[Corollary 5.15.]{compositionpaper});
\item[$(ii)$] \hyperlink{om1pp}{$(\omega_{1'})$}, i.e. $\omega$ is equivalent w.r.t. \hyperlink{sim}{$\sim$} to a sub-additive weight, see \cite[Theorems 3,4]{characterizationstabilitypaper} and \cite[Chapter 6]{compositionpaper}.
\end{itemize}

\subsection{Weight matrices in the sense of Beaugendre, Schmets and Valdivia}\label{section63}
Beaugendre in \cite{Beaugendre01} and Schmets and Valdivia in \cite{intersectionpaperextension} have considered weight matrices in the following sense: Let $\Phi:[0,+\infty)\rightarrow\RR$ be a convex and increasing function with $\lim_{t\rightarrow\infty}\frac{\Phi(t)}{t}=+\infty$ and $\Phi(0)=0$ (w.l.o.g. - replace $\Phi$ by $\Psi(t):=\Phi(t)-\Phi(0)$, see \cite[Definition 16.]{intersectionpaperextension}). We introduce the following weight matrix
$$\mathcal{M}^{\Phi}:=\{(p!m^{\Phi}_{ap})_{p\in\NN}: a>0\}\hspace{40pt}m^{\Phi}_{ap}:=\exp(\Phi(ap)).$$
In the literature the Beurling-case $\mathcal{E}_{(\mathcal{M}^{\Phi})}$ was considered. We summarize some properties:

\begin{itemize}
\item[$(i)$] $\mathcal{M}^{\Phi}$ is \hyperlink{Msc}{$(\mathcal{M}_{\on{sc}})$} and \hyperlink{B-Comega}{$(\mathcal{M}_{(\mathcal{C}^{\omega})})$} holds.
\item[$(ii)$] \hyperlink{R-L}{$(\mathcal{M}_{\{\text{L}\}})$} and \hyperlink{B-L}{$(\mathcal{M}_{(\text{L})})$} both are satisfied, compare this with \cite[Lemma 5.9. (5.10)]{compositionpaper} where condition \hyperlink{om1}{$(\omega_1)$} is needed. As shown in \cite[Lemma 17]{intersectionpaperextension} we get both
    $$\forall\;a>0\;\forall\;h>0\;\exists\;b>0\;(b>a)\;\exists\;D>0\;\forall\;p\in\NN_{>0}:\;\;\log(h)-\frac{\log(D)}{p}\le\frac{1}{p}\left(\Phi(bp)-\Phi(ap)\right)$$
    and
    $$\forall\;b>0\;\forall\;h>0\;\exists\;a>0\;(a<b)\;\exists\;D>0\;\forall\;p\in\NN_{>0}:\;\;\log(h)-\frac{\log(D)}{p}\le\frac{1}{p}\left(\Phi(bp)-\Phi(ap)\right),$$
    since convexity of $\Phi$ yields
    \begin{equation}\label{Phiconvexity}
    \forall\;a,b>0,\;b>a:\;\frac{\Phi(bp)-\Phi(ap)}{p(b-a)}\ge\frac{\Phi(bp)}{pb}\Leftrightarrow\frac{\Phi(bp)-\Phi(ap)}{p}\ge\frac{\Phi(bp)}{pb}(b-a)\rightarrow\infty
    \end{equation}
    as $p\rightarrow\infty$.
\item[$(iii)$] \eqref{Phiconvexity} implies also that all sequences are pairwise not equivalent. If $(m^{\Phi}_{ap})_p\hyperlink{approx}{\approx}(m^{\Phi}_{bp})_p$ for all $a,b>0$, then we would get
$$\forall\;a>0\;\forall\;b>0\;\exists\;C\ge 1\;\forall\;p\in\NN:\;\;m^{\Phi}_{bp}\le C^p m^{\Phi}_{ap}\Leftrightarrow\frac{1}{p}\left(\Phi(bp)-\Phi(ap)\right)\le\log(C),$$
but the left hand side tends to infinity as $p\rightarrow\infty$ whenever $b>a$. So $\mathcal{M}^{\Phi}$ has both \hyperlink{R-strict}{$(\mathcal{M}_{\{\text{strict}\}})$} and \hyperlink{B-strict}{$(\mathcal{M}_{(\text{strict})})$}.
\item[$(iv)$] $\mathcal{M}^{\Phi}$ has \hyperlink{R-mg}{$(\mathcal{M}_{\{\text{mg}\}})$} and \hyperlink{B-mg}{$(\mathcal{M}_{(\text{mg})})$}. By convexity of $\Phi$ we get $\Phi(ap+aq)\le\frac{1}{2}\Phi(2ap)+\frac{1}{2}\Phi(2aq)\le\Phi(2ap)+\Phi(2aq)$ for all $a>0$ and $p,q\in\NN$ and so
    $$M^{\Phi}_{a(p+q)}\le M^{\Phi}_{bp}\cdot M^{\Phi}_{bq}\Leftrightarrow\Phi(a(p+q))\le\Phi(bp)+\Phi(bq)$$
    holds with $b=2a$.
\item[$(v)$] \hyperlink{R-FdB}{$(\mathcal{M}_{\{\text{FdB}\}})$} and \hyperlink{B-FdB}{$(\mathcal{M}_{(\text{FdB})})$} both are satisfied. This is clear since each $(m^{\Phi}_{ap})_p$ is log-convex, see e.g. \cite[2.2. Lemma $(1)$]{compositionpaper}.
\end{itemize}

Thus also for $\mathcal{E}_{[\mathcal{M}^{\Phi}]}$ the exponential laws in Theorem \ref{exponentiallaw} and the consequences in Lemma \ref{expconsequence1} are valid. Moreover the characterizing results \cite[Theorems 5,6]{characterizationstabilitypaper} and all further generalizations of the results from \cite{KMRplot} hold.\vspace{6pt}

As special case one may consider $\Phi=\varphi^{*}_{\omega}$ for $\omega\in\hyperlink{omset1}{\mathcal{W}}$. Then on the one hand one has the matrix $\mathcal{M}^{\Phi}$ as defined before, on the other hand the weight matrix $\Omega:=\{(\Omega^a_p)_p:=\exp(1/a\varphi^{*}_{\omega}(ap)): a>0\}$ as the approach in \cite{BraunMeiseTaylor90}. By definition we have
$$m^{\Phi}_{ap}:=\exp(\Phi(ap))=\exp(1/a\varphi^{*}_{\omega}(ap))^a=(\Omega^a_p)^a.$$
As we have already pointed out the weights $\omega_s:=\max\{0,\log(t)^s\}$, $s>1$, generate an infinite non-constant weight matrix. We denote the associated matrices by $\mathcal{M}_s^{\Phi}$ and $\Omega_s$ and prove:

\begin{lemma}
For any $s>1$ the matrices $\mathcal{M}^{\Phi}_s$ and $\Omega_s$ are equivalent w.r.t. both \hyperlink{Mroumapprox}{$\{\approx\}$} and \hyperlink{Mbeurapprox}{$(\approx)$}.
\end{lemma}

\demo{Proof}
Let $s>1$ be arbitrary but fixed. For $t\ge 0$ we get $\varphi_{\omega_s}(t)=\omega_s(\exp(t))=(\log(\exp(t)))^s=t^s$, hence $\varphi^{*}_{\omega_s}(x)=\sup\{xy-y^s:y\ge 0\}=:\sup\{f_{x,s}(y): y\ge 0\}$ for all $x\ge 0$. A straightforward computation shows
$$\varphi^{*}_{\omega_s}(x)=f_{x,s}\left(\left(\frac{x}{s}\right)^{\frac{1}{s-1}}\right)=x\left(\frac{x}{s}\right)^{\frac{1}{s-1}}-\left(\frac{x}{s}\right)^{\frac{s}{s-1}}=x^{\frac{s}{s-1}}\left(\frac{1}{s^{\frac{1}{s-1}}}-\frac{1}{s^{\frac{s}{s-1}}}\right)=:x^{\frac{s}{s-1}}R(s)$$ and so
\begin{equation}\label{examplematrixsequence}
\Omega^l_p=\exp\left(l^{1/(s-1)}p^{s/(s-1)}R(s)\right)\hspace{30pt}m^{\Phi}_{lp}=\exp(l^{s/(s-1)}p^{s/(s-1)}R(s)).
\end{equation}
The case $s=2$ gives $\Omega^l_p=(\exp(lR(2)))^{p^2}=(\exp(l/4))^{p^2}$ and $m^{\Phi}_{lp}=(\exp(l^2/4))^{p^2}$.

$\Omega_s\hyperlink{Mroumpreceq}{\{\preceq\}}\mathcal{M}^{\Phi}_s$. Let $l\in\NN_{>0}$ (large) and get $\Omega^l_p\le(\Omega^l_p)^l\le(\Omega^l_p)^lp!=p!m^{\Phi}_{lp}$ for each $p\in\NN$ since $\Omega^l_p\ge 1$ for each $l>0$, $p\in\NN$.

$\mathcal{M}^{\Phi}_s\hyperlink{Mroumpreceq}{\{\preceq\}}\Omega_s$. Let $l>0$, then we have to find $n>l>0$ and $C\ge 1$ such that for all $p\in\NN$ we get $p!m^{\Phi}_{lp}\le C^p\Omega^n_p\Leftrightarrow p!\exp(l^{s/(s-1)}p^{s/(s-1)}R(s))\le C^p\exp(n^{1/(s-1)}p^{s/(s-1)}R(s))$. So the choice $n=2^{s-1}l^s$ is sufficient and analogously $\mathcal{M}^{\Phi}_s\hyperlink{Mbeurpreceq}{(\preceq)}\Omega_s$ holds, too.

$\Omega_s\hyperlink{Mbeurpreceq}{(\preceq)}\mathcal{M}^{\Phi}_s$. For each $l>0$ (small) there exists $C\ge 1$ and $n>0$ such that for all $p\in\NN$ we get $\Omega^n_p\le C^pp!m^{\Phi}_{lp}\Leftrightarrow\exp(n^{1/(s-1)}p^{s/(s-1)}R(s))\le C^p p!\exp(l^{s/(s-1)}p^{s/(s-1)}R(s))$, so the choice $n=l^s$ is sufficient.
\qed\enddemo

If $\omega\in\hyperlink{omset1}{\mathcal{W}}$, then $\Omega$ has always both \hyperlink{R-mg}{$(\mathcal{M}_{\{\text{mg}\}})$} and \hyperlink{B-mg}{$(\mathcal{M}_{(\text{mg})})$}. But $\Omega^l\hyperlink{approx}{\approx}\Omega^n$ for all $l,n>0$ holds if and only if \hyperlink{mg}{$(\text{mg})$} for some/each $\Omega^l$ and if and only if \hyperlink{om6}{$(\omega_6)$} for $\omega$, see \cite[Chapter 5]{compositionpaper}.

For $\mathcal{M}^{\Phi}$ this is not true any more. As we have already seen the sequences in $\mathcal{M}^{\Phi}$ are always pairwise not equivalent.

On the other hand, since $(m^{\Phi}_{ap})_{p\in\NN}$ is log-convex, \hyperlink{mg}{$(\text{mg})$} holds for this sequence if and only if $m^{\Phi}_{a2p}\le C^{2p}(m^{\Phi}_{ap})^2\Leftrightarrow\frac{1}{2p}\Phi(2ap)-\frac{1}{p}\Phi(ap)\le\log(C)$ for a constant $C\ge 1$ and all $p\in\NN$, see \cite[Theorem 1, $(3)\Rightarrow(2)$]{matsumoto}. So if $\Phi$ satisfies
\begin{equation}\label{Phimg}
\exists\;D\ge 1\;\forall\;t\ge 0:\;\Phi(2t)\le 2\Phi(t)+Dt,
\end{equation}
then each $(m^{\Phi}_{ap})_{p\in\NN}$ has \hyperlink{mg}{$(\text{mg})$}. In \cite{Beaugendre01} a weight with \eqref{Phimg} is called a {\itshape weight of moderate growth.} \eqref{Phimg} is valid for $\Phi=\varphi^{*}_{\omega}$ if and only if $\omega\in\hyperlink{omset1}{\mathcal{W}}$ has \hyperlink{om6}{$(\omega_6)$}. This holds by the proof of $(5.11.)$ in \cite[Lemma 5.9.]{compositionpaper} and by applying the conjugate operator to \eqref{Phimg} (note that $\varphi^{**}_{\omega}=\varphi_{\omega}$).\vspace{6pt}

Finally consider $\Phi(t):=t\log(t)$ for $t\ge 1$ and $\Phi(t):=0$ for $0\le t<1$. Each $(m^{\Phi}_{ap})_{p\in\NN}$ has \hyperlink{mg}{$(\text{mg})$}, since $\frac{1}{2p}\Phi(2ap)-\frac{1}{p}\Phi(ap)=a\log(2)$. More precisely Stirling's formula and $m^{\Phi}_{ap}=\exp(\Phi(ap))=(ap)^{ap}$ show that this yields the Gevrey-matrix $\mathcal{G}$ and which should be compared with \cite[5.19.]{compositionpaper}.

\bibliographystyle{plain}
\bibliography{Bibliography}
\end{document}